\documentclass[12pt]{article}
\usepackage{amsmath}
\usepackage{amsthm}
\usepackage{caption}
\usepackage{amssymb,enumerate}
\usepackage{latexsym}
\usepackage{cite}
\usepackage{indentfirst}
\usepackage{xcolor}
\usepackage{enumitem}
\newif\ifshowrevisions
\showrevisionsfalse
\ifshowrevisions
  \newenvironment{revisedblock}{\color{blue}}{}
  \newcommand{\revnote}[1]{\par\noindent{\color{red}\small\textbf{[Revision note: #1]}}\par}
\else
  \newenvironment{revisedblock}{}{}
  \newcommand{\revnote}[1]{}
\fi
\usepackage{amsmath}
\usepackage{verbatim}
\usepackage{graphicx} 
\usepackage{epstopdf}
\usepackage{float} 
\usepackage{subfigure}
\usepackage{mathrsfs}
\usepackage{epstopdf}
\newcommand{\R}{\mathbb{R}}

\newtheorem{theorem}{Theorem}[section]
\newtheorem{lemma}[theorem]{Lemma}
\newtheorem{proposition}[theorem]{Proposition}

\newtheorem{corollary}[theorem]{Corollary}

\usepackage[colorlinks=true,
linkcolor=black, 
anchorcolor=black,  
citecolor=black,  
urlcolor=blue 
]{hyperref}

\newcommand{\DOI}[1]{doi:\href{https://doi.org/#1}{#1 }}

\newtheorem{remark}[theorem]{Remark}
\newcommand{\be}{\begin{equation}}
	\newcommand{\ee}{\end{equation}}

\makeatletter

\newcommand{\Rmnum}[1]
{expandafter\@slowromancap\romannumeral#1@}
\makeatother
\numberwithin{figure}{section}
\allowdisplaybreaks[4]
\usepackage{amssymb}
\usepackage{amsmath}
\begin{document}
\title{{\Large\bf  $L^2$-contraction and Asymptotic Stability of Monotone Shock Waves for Cubic and Quartic KdV-Burgers Equations}}
\author{Yi Du,\ \ Tao Yan, \ \ Yushu Zeng, \ \ Shihui Zhu\thanks{Corresponding author \textit{E-mail address}: shihuizhumath@sicnu.edu.cn}\,\,\,\,\,\,\small\\
{\small \textit{School of Mathematical Sciences, Sichuan Normal University}}\\
{\small \textit{Chengdu 610066, PR China}}}
\date{}
 \maketitle

\begin{center}
\begin{minipage}{6.0in}
{\textbf{Abstract:} We investigate the nonlinear stability of monotone traveling shock waves arising in the generalized Korteweg–de Vries–Burgers equation with a power-law flux \(f(\phi)=\phi^p/p\), focusing on the cubic and quartic cases \(p=3\) and \(p=4\), respectively, in the presence of both dissipative and dispersive effects.
Since shocks connecting distinct end states destroy the integrable structure of the equation, we introduce a slowly varying reference connection to establish global well-posedness for the Cauchy problem in an affine Sobolev space. A key ingredient in the stability analysis is to establish a new two-sided comparison estimate for the derivative of a monotone shock-wave profile. After a monotone coordinate transformation, this estimate allows the viscous dissipation to be compared with the natural vanishing weight arising in a shifted
Poincar{\'e} inequality. Together with an $L^\infty$ bootstrap controlling the higher-order relative-flux terms, this yields a global $L^2$-contraction property modulo a time-dependent shift for sufficiently small $H^1$ perturbations in the admissible parameter regimes. We  obtain exponential decay of the shock 
profile and show asymptotic convergence to the shifted shock in $L^q(\R)$, $2<q\leq \infty$, with vanishing shift velocity.}  

\end{minipage}
\end{center}
\vspace{0.3cm}

\renewcommand{\theequation}{\thesection.\arabic{equation}}
\setcounter{equation}{0}

	\section{Introduction}

The one-dimensional generalized Korteweg-de Vries-Burgers equation (gKdV-Burgers equation) is given by 
\begin{align}\label{eq}
    \begin{cases}
	\phi_t+\left(f(\phi)\right)_x=\epsilon \phi_{xx}-\delta \phi_{xxx},\qquad t\in \R_+,\quad x\in \R,\\
	\phi(0,x)=\varphi(x),
	\end{cases}
\end{align}
where 
 $\epsilon>0$ and $\delta>0$ are the coefficients of dissipation and dispersion. $f(\phi)=\frac{\phi^p}{p}$ is the nonlinear flux.  Indeed, the classical quadratic case $p=2$ has its origins in the long-wave theories of Boussinesq \cite{Boussinesq}, Korteweg–de Vries \cite{DJkorteweg1895}, and Su and Gardner \cite{SuGardner1969}. 
The gKdV-Burgers equation serves as a canonical model to describe the competition among nonlinear steepening, viscous dissipation and dispersion. 
This interplay can give rise to coherent shock structures  whose shape and long-time behavior depend sensitively on 
the relative strength of dissipation and dispersion (see \cite{Chen2000}). 

In many weakly nonlinear asymptotic reductions, the power of the nonlinear flux is determined by the first nonvanishing nonlinear term when lower-order nonlinear coefficients vanish. In generic regimes, the leading nonlinear coefficient is nonzero, yielding 
the conventional term $\phi\phi_x$, which corresponds to $p=2$ in the present notation. At a critical composition, 
however, the quadratic nonlinear coefficient may vanish. Thus, the standard KdV reduction loses its leading nonlinear balance and the asymptotic expansion  
must be carried to higher order. This leads naturally to the modified KdV equation with cubic nonlinearity \cite{p3Demiray1999,p3Verheest2019,p3PUXI2021}, 
and, in the presence of dissipation, to the modified KdV-Burgers equation \cite{p3Demiray1999,p3Sarma2009}. At a supercritical composition, even the cubic coefficient may vanish 
together with the quadratic one, leaving a quartic nonlinearity as the leading nonlinear contribution \cite{p4Verheest2016}. 

Motivated by these higher-order asymptotic regimes, we focus on the cases $p\in\{3,4\}$, for which the structure of traveling waves and the mechanisms underlying nonlinear stability may differ substantially from those in the classical quadratic case. Moreover, the two cases exhibit distinct symmetries under sign reversal. In the case $p=3$, changing the sign of both the solution and the initial datum preserves the governing equation. More precisely, a solution $\phi(t,x)$ with initial datum $\varphi(x)$ gives rise to another solution $-\phi(t,x)$ associated with $-\varphi(x)$, while the coefficients $\varepsilon$ and $\delta$ remain fixed. In contrast, this symmetry is lost when $p=4$, since the nonlinear term $\phi^{3}\phi_x$ remains unchanged under $\phi\mapsto-\phi$, whereas the time derivative, dissipative term, and dispersive term all change sign. Thus, although $p=3$ and $p=4$ represent the first two higher-order nonlinearities beyond the quadratic case, their associated equations possess distinct symmetry properties, which may lead to differences in the analysis of traveling waves and their stability.

Before turning to the construction and stability of shock waves, it is necessary to establish a 
global solution framework for the corresponding Cauchy problem. Since the shock-type solutions 
considered here connect two distinct nonzero end states,
\begin{align*}
	\lim_{x\to -\infty} \phi(t,x)=\phi_{-\infty},\qquad \lim_{x\to +\infty} \phi(t,x)=\phi_{+\infty},
\end{align*}
the function $\phi(t,\cdot)$ itself does not belong to $H^1(\R)$. It is therefore natural to seek solutions in the affine Sobolev space
$\bar{\phi}+H^1(\mathbb{R})$, where $\bar{\phi}$
is a fixed smooth function connecting the
two end states. To specify this space, let $\chi\in C^\infty(\R)$ 
be a nondecreasing cutoff function satisfying 
\begin{align*}
	\chi(x)=0 \quad \text{for} \quad x\leq -1,\qquad \chi(x)=1 \quad \text{for}\quad  x\geq 1.
\end{align*}
We fix a smooth reference connection $\bar{\phi}$ by 
\begin{align}\label{baru}
	\bar{\phi}:=\phi_{-\infty}+(\phi_{+\infty}-\phi_{-\infty})\chi(x).
\end{align}
Then $\bar{\phi}$ is in $C^\infty(\R)\cap L^\infty(R)$, and all its derivatives of positive order are compactly supported. 
We emphasize that $\bar{\phi}$ serves only to specify the affine Sobolev space and is not assumed to be a traveling-wave profile. 
Classical and modified KdV-Burgers equations have been studied extensively in 
the zero-background or periodic settings, including low-regularity and global results \cite{Molinetlowregularity2002,Zhang2009}, while generalized nonlinearities 
have been considered under suitable growth assumptions \cite{Gallego2019}. For distinct nonzero end states, a related quadratic theory has recently been developed 
in \cite{chen2026}. In the higher-order cases considered here, passing to a perturbation
variable alone does not suffice to close the global energy estimate,  
because superlinear powers of the perturbation are produced by the nonlinear flux.

To overcome this difficulty, we exploit the freedom in the choice of the reference connection and establish global 
well-posedness for the gKdV–Burgers equation with $p=3,4$ around two distinct end states. Combining the resulting $L^2$
control with a modified $H^1$ energy estimate yields global continuation of the local solution.

\begin{theorem}
	Let $p\in\{3,4\}$ and $\epsilon,\delta>0$. Denote $\bar{\phi}$ as the fixed reference connection defined by \eqref{baru}. If 
	$\varphi-\bar{\phi}\in H^1(\R)$, then the Cauchy problem \eqref{eq} admits a unique global solution satisfying
	\begin{align*}
		\phi-\bar{\phi}\in C\bigl([0,\infty);H^1(\R)\bigr)
\cap L^2_{\mathrm{loc}}\bigl([0,\infty);H^2(\R)\bigr).
	\end{align*}
\end{theorem}

\begin{remark}
	The above global well-posedness result provides an essential foundation for the nonlinear stability analysis of shock 
	waves developed below. Indeed, a global-in-time stability statement is meaningful only when the perturbed Cauchy problem 
	admits a global solution in the natural shock phase space. The $H^1$ bounds for the perturbation, together with the one-dimensional Sobolev embedding and the boundedness of the reference profile, yield $L^\infty$ bounds for the solution on finite time intervals. These bounds play an important role in controlling the higher-order relative-flux terms arising in the cases $p=3$ and $p=4$. Thus, global well-posedness supplies both the global solution framework and the regularity needed for the subsequent shifted contraction estimates.
\end{remark}


To investigate such shock waves, we seek
traveling-wave solutions of (1) of the form
\begin{align*}
\phi(t,x)=\widetilde{\phi}(\xi),\qquad \xi=x-\sigma t,
\end{align*}
where $\sigma$ denotes the wave speed. We consider shock profiles connecting two distinct constant end states,
\begin{align}\label{boundcondition}
\lim_{\xi\to-\infty}\widetilde{\phi}(\xi)=\phi_{-\infty},
\qquad
\lim_{\xi\to+\infty}\widetilde{\phi}(\xi)=\phi_{+\infty},
\qquad
\phi_{-\infty}>\phi_{+\infty}.
\end{align}
Substitution into \eqref{eq} gives the traveling-wave equation
\begin{align}\label{eq:shockwave}
-\sigma\widetilde{\phi}'
+\widetilde{\phi}^{p-1}\widetilde{\phi}'
=
\epsilon\widetilde{\phi}''
-\delta\widetilde{\phi}'''.
\end{align}
After one integration in $\xi$, this becomes
\begin{align*}
\delta \widetilde{\phi}''
-\epsilon\widetilde{\phi}'
+    \frac{\phi^p-\phi_{+\infty}^p}{p}
    -\sigma(\phi-\phi_{+\infty})=0,
\end{align*}
and the end conditions determine the Rankine-Hugoniot speed
\begin{align}\label{wavespeed}
\sigma
=
\frac{\phi_{-\infty}^p-\phi_{+\infty}^p}
{p(\phi_{-\infty}-\phi_{+\infty})}.
\end{align}

Among viscous–dispersive shock waves, monotone profiles form a particularly important class, as they provide nonoscillatory transitions between the two far-field states and preserve the ordering of these states throughout the shock layer. These profiles constitute an important class of solutions in the study of nonlinear stability for viscous–dispersive shock waves.

The existence of monotone traveling-wave profiles for generalized KdV-Burgers equations has been studied in 
\cite{Chen2000}. For the power fluxes considered here, 
the known monotonicity criterion can be stated in the following form.

\begin{proposition}( \cite{Chen2000})\label{prop:monotone-profile}
Let $p\in\{3,4\}$ and $\epsilon,\delta>0$. Assume
$\phi_{-\infty}>\phi_{+\infty}$ and 
\begin{align}\label{monotone-profile-condition}
\epsilon^2
\geq
4\delta\bigl(\phi_{-\infty}^{p-1}-\sigma\bigr),
\end{align}
with $\sigma$ given by \eqref{wavespeed}. In the cubic case $p=3$, assume in addition that
$\phi_{+\infty}>0$. Under these conditions, one can find a traveling-wave
profile $\widetilde{\phi}$ satisfying
\eqref{boundcondition} and
$\widetilde{\phi}'(\xi)<0$ for every $\xi\in\mathbb{R}$.
The profile is unique up to translation.
\end{proposition}

For the cubic and quartic fluxes, condition \eqref{monotone-profile-condition} becomes, respectively,
\begin{align*}
\epsilon^2
&\geq
\frac{4\delta}{3}
(\phi_{-\infty}-\phi_{+\infty})
(2\phi_{-\infty}+\phi_{+\infty}),
&& p=3,\\
\epsilon^2
&\geq
\delta
(\phi_{-\infty}-\phi_{+\infty})
\bigl(
3\phi_{-\infty}^2
+2\phi_{-\infty}\phi_{+\infty}
+\phi_{+\infty}^2
\bigr),
&& p=4.
\end{align*}
Notice that the parameter assumptions imposed in the contraction theorems below imply \eqref{monotone-profile-condition}. 
Hence the monotone shock profiles required in our stability analysis are available throughout the parameter regimes considered 
in this paper.




Even after restricting to monotone shocks, there is a fundamental obstruction to stability with respect to a fixed profile. 
Traveling waves are translation invariant. In particular, a solution may differ from the reference shock only by a spatial displacement and yet remain 
a finite distance from the unshifted profile for all time. Indeed, even the exact translated wave 
$\phi(t,x)=\widetilde{\phi}(x-\sigma t-a)$
does not, in general, satisfy 
\begin{align*}
	\|\phi(t,\cdot)-\widetilde{\phi}(\cdot-\sigma t)\|_{L^2(\R)}\to 0,\qquad t\to \infty
\end{align*}
when $a\neq 0$. At the linearized level, this degeneracy is represented by the neutral translational mode 
$\widetilde{\phi}'$. Therefore the natural stability notion is orbital stability with respect to the family 
of translates rather than convergence to one fixed representative. This leads to the introduction of a time-dependent shift 
$Y(t)$ and to the comparison of $\phi(t,x)$ with $\widetilde{\phi}(x-\sigma t-Y(t))$. The shift is 
thus forced by the geometry of the traveling-wave manifold rather than introduced merely as a 
technical device (see \cite{kang2017, chen20262}).


This translational difficulty is already present in the classical quadratic KdV-Burgers equation. Traveling shock profiles for 
the quadratic equation were studied by Bona and Schonbek \cite{Bona1985}, while Pego established asymptotic stability in a 
perturbative regime \cite{Pego1985}. Generalized power nonlinearities were subsequently investigated from a spectral 
viewpoint in \cite{Pego1993}. Zhang et al. \cite{zhang2019} established the asymptotic stability of monotone decreasing kink profiles for gKdV–Burgers equations 
under smallness and antiderivative assumptions on the perturbation. For scalar viscous conservation laws, Kang and Vasseur \cite{kang2017} developed a shift-based $L^2$ 
contraction method that treats the neutral translational mode at the nonlinear level. In a subsequent study, Chen et al. \cite{chen20262} applied this framework
to the quadratic KdV-Burgers equation, demonstrating that monotone shocks
satisfy an $L^2$ contraction property under a suitable time-dependent shift. Eun et al. \cite{Eun2026} analyzed the long-time stability of degenerate shock-rarefaction composite 
waves for the cubic modified KdV–Burgers equation via the $a$-contraction method with shifts. 
 The quadratic problem therefore provides the natural benchmark for the present work. Our purpose is to determine to what 
extent this contraction mechanism survives when the quadratic flux is replaced by the genuinely higher-order fluxes $p\in \{3,4\}$. 

\begin{theorem}[\(L^2\)-contraction for \(p = 3\)]
	\label{thm:main}
	Let $\epsilon,\delta>0$ and
	$
	\phi_{-\infty}>\phi_{+\infty}>0$ satisfy 
	$\phi_{-\infty}<\dfrac{8(3\sqrt{2}+1)}{17}\phi_{+\infty},
	$
	and assume that
	\[
	0<
	\frac{\delta(\phi_{-\infty}^{2}-\sigma)}{\epsilon^{2}}
	\le \frac14.
	\]
	Denote by $\widetilde{\phi}$ the associated monotone
	decreasing shock profile joining $\phi_{-\infty}$
	to $\phi_{+\infty}$, and by $\phi$ the global solution
	to \eqref{eq} emanating from the initial datum $\varphi$.
	
	Fix $\gamma$ satisfying
	\(
	0<\gamma<8(\sqrt{2}-1)\phi_{+\infty}
	+(4\sqrt{2}-7)\phi_{-\infty}.
	\)
	There is a constant $\varepsilon_0>0$ such that,
	whenever
	\(
	\|\varphi-\widetilde{\phi}\|_{H^1(\mathbb{R})}
	<\varepsilon_0,
	\)
	the time-dependent shift $Y(t)$ is specified by
	the initial condition $Y(0)=0$ together with
	\[
	\dot{Y}(t)
	=
	\frac{c_{1}}{\phi_{-\infty}-\phi_{+\infty}}
	\int_{\mathbb{R}}
	\left[
	\phi(t,x)
	-\widetilde{\phi}(x-\sigma t-Y(t))
	\right]
	\widetilde{\phi}'(x-\sigma t-Y(t))\,dx,
	\]
	where
	$
	c_{1}
	=
	-\frac{
		\gamma
		+3\phi_{-\infty}
	}{3}. 
	$
	Moreover, for any $t\in \mathbb{R}_+$,
	\[
	\begin{aligned}
		\left\|
		\phi(t,\cdot)
		-\widetilde{\phi}(\cdot-\sigma t-Y(t))
		\right\|_{L^{2}(\mathbb{R})}^{2}
		&
		+2A_{1}^{*}\epsilon
		\int_{0}^{t}
		\left\|
		\partial_x
		\left[
		\phi(s,\cdot)
		-\widetilde{\phi}(\cdot-\sigma s-Y(s))
		\right]
		\right\|_{L^{2}(\mathbb{R})}^{2}\,ds
		\\
		&\le
		\|\varphi-\widetilde{\phi}\|_{L^{2}(\mathbb{R})}^{2},
	\end{aligned}
	\]
	where
	$
	A_{1}^{*}
	=
	\frac{8(\sqrt{2}-1)\phi_{+\infty}+(4\sqrt{2}-7)\phi_{-\infty}-\gamma}{6\lambda_{2}}. 
	$
	
\end{theorem}

\begin{remark}
The passage from $p=2$ to $p=3,4$ is not a routine modification of the exponent. The first new difficulty appears 
directly in the relative flux. Writing $w=\phi(x-\sigma t+Y(t))-\widetilde{\phi}$, the nonlinear term in the $L^2$ energy identity contains 
$f(\widetilde{\phi}+w)-f(\widetilde{\phi})-f'(\widetilde{\phi})w$. For the quadratic flux, this expression is exactly 
quadratic in $w$, which matches the quadratic structure of the viscous dissipation. For $p=3$ and $p=4$, the nonlinear terms are 
$\widetilde{\phi}w^2+\frac13w^3$ and $\frac32\widetilde{\phi}^2w^2+\widetilde{\phi}w^3+\frac14w^4$. Thus, genuinely 
cubic and quartic contributions enter the contraction identity and destroy the direct quadratic balance available for $p=2$. Their control requires an additional $L^\infty$ bootstrap for the perturbation together with the global $H^1$ 
information developed in this paper.
\end{remark}

A second, independent difficulty concerns the profile weight generated by the viscous term. After the monotone change of 
variables $z=\widetilde{\phi}(\xi)$, the dissipation takes the form $\int (-\epsilon\widetilde{\phi}')|v_z|^2dz.$ 
To combine this term with the shifted Poincar{\'e} inequality, one needs a quantitative two-sided comparison between 
$-\epsilon\widetilde{\phi}'$ and the natural vanishing weight $(\phi_{-\infty}-\widetilde{\phi})(\widetilde{\phi}-\phi_{+\infty})$. 
The auxiliary comparison structure used in the quadratic case does not retain the required global behavior for higher powers. This difficulty 
is particularly significant in the quartic problem. We overcome it by introducing a new comparison function for the profile equation 
and proving a general two-sided estimate.

Having established the $L^2$-contraction property for the cubic case $p=3$, we next turn to the quartic 
case $p=4$. Although the overall shifted energy framework remains applicable, the quartic nonlinearity 
introduces additional higher-order terms and requires a corresponding refinement of the estimates used in 
the cubic case. We therefore investigate whether the contraction mechanism can be extended to this
higher-order nonlinearity. The resulting $L^2$-contraction estimate for $p=4$ is stated as follows.

		\begin{theorem}[$L^2$-contraction for $p=4$]\label{the:p'=4}
		Let $\epsilon,\delta>0$ and $\phi_{-\infty}>\phi_{+\infty}$. Assume that
		\[
		0<\frac{\delta(\phi_{-\infty}^3-\sigma)}{\epsilon^2}\le\frac14.
		\]
	Consider the solution $\phi$ to \eqref{eq} arising
	from the initial datum $\varphi$, together with
	the associated strictly decreasing shock profile
	$\widetilde{\phi}$ joining the end states
	$\phi_{-\infty}$ and $\phi_{+\infty}$.
	
	Define
	\(
	U_*=\max\{|\phi_{-\infty}|,|\phi_{+\infty}|\},
	\)
	and impose the condition
	\(
	\kappa_1m_4>3U_*^2.
	\)
	For any $\gamma$ satisfying
	\(
	0<\gamma<
	\sqrt{2(\kappa_1m_4-U_*^2)}-2U_*,
	\)
	there is a constant $\varepsilon_0>0$ such that,
	whenever
	\(
	\|\varphi-\widetilde{\phi}\|_{H^1(\mathbb{R})}
	<\varepsilon_0,
	\)
	the time-dependent shift $Y(t)$ is prescribed
	by the initial condition $Y(0)=0$ and
		\[
		\dot Y(t)
		=
		\frac{c_2}
		{\phi_{-\infty}-\phi_{+\infty}}
		\int_{\mathbb R}
		\bigl[\phi(t,x)-\widetilde \phi(x-\sigma t-Y(t))\bigr]
		\widetilde \phi'(x-\sigma t-Y(t))\,dx,
		\]
		where $c_2=-(\frac{\gamma^2}{4}+\gamma U_*+\frac32U_*^2)$.	Then, for any $t\in\mathbb{R}_+$,
		\begin{equation*}
			\begin{aligned}
				\|\phi(t,\cdot)-\widetilde \phi(\cdot-\sigma t-Y(t))\|_{L^2(\R)}^2
				&
				+2A_2^*\epsilon
				\int_0^t
				\|\partial_x[\phi(s,\cdot)
				-\widetilde \phi(\cdot-\sigma s-Y(s))]\|_{L^2(\R)}^2\,ds\\
				&\le
				\|\varphi(x)-\widetilde \phi(x)\|_{L^2(\R)}^2,
			\end{aligned}
		\end{equation*}
		where
		$
		A_2^*
		=
		\frac{
			2\lambda_3-
			\left(\frac{\gamma^2}{4}+\gamma U_*+\frac32U_*^2\right)
		}{2\lambda_4}>0.
		$
	\end{theorem}

\begin{remark}
	\begin{enumerate}[label=(\roman*), ref=(\roman*), leftmargin=*, labelindent=0pt, itemsep=0pt, topsep=0.25em]
\item 	The parameters in the assumptions of Theorem \ref{the:p'=4} are defined explicitly as follows:  
\begin{align*}
	m_4
	&:= \min_{z\in[\phi_{+\infty},\phi_{-\infty}]}H(z),
	\qquad
	M_4
	:= \max_{z\in[\phi_{+\infty},\phi_{-\infty}]}H(z),
	\qquad
	\lambda_3
	= \frac{\kappa_1 m_4}{4},
	\qquad
	\lambda_4
	= \frac{\kappa_2 M_4}{4},\\
	\kappa_1
	&= \frac{2}{
		1+\sqrt{1+\dfrac{4\delta\bigl(\sigma-\phi_{+\infty}^{p-1}\bigr)}{\epsilon^2}}
	},
	\qquad
	\kappa_2
	= \frac{2}{
		1+\sqrt{1-\dfrac{4\delta\bigl(\phi_{-\infty}^{p-1}-\sigma\bigr)}{\epsilon^2}}
	},
\end{align*}
	where $	H(z)
	= z^2+(\phi_{-\infty}+\phi_{+\infty})z
	+\phi_{-\infty}^2+\phi_{-\infty}\phi_{+\infty}+\phi_{+\infty}^2$. 
	
\item 	The condition \(\kappa_1 m_4 > 3U_*^2\) excludes shocks whose end states cross zero, but it can be satisfied when both end states are positive or both are negative, subject to the stated parameter restrictions.

\item We illustrate that the admissible parameter set is nonempty. Let
$$
\phi_{-\infty}=1.2,\qquad
\phi_{+\infty}=1,\qquad
\epsilon=1,\qquad
\delta=0.01.
$$
Then the Rankine–Hugoniot speed is \(\sigma=1.342\), and
$$
0<
\frac{\delta(\phi_{-\infty}^{3}-\sigma)}{\epsilon^2}
=0.00386<\frac14.
$$
Moreover, direct calculations yield $
\kappa_1\approx0.996603$, $
m_4=6.84$, $ U_*=1.2,
$
so that
$$
\kappa_1m_4\approx6.816766
>3U_*^2=4.32.
$$
Consequently,
$
\sqrt{2(\kappa_1m_4-U_*^2)}-2U_*
\approx0.879258>0$. 
Choosing \(\gamma=0.2\), we further obtain
$$
2\lambda_3-
\left(\frac{\gamma^2}{4}+\gamma U_*+\frac32U_*^2\right)
\approx0.998383>0.
$$
Thus, the parameter conditions required in Theorem \ref{the:p'=4} can be satisfied simultaneously, showing that the admissible parameter set is nonempty.

\end{enumerate}
\end{remark}

The results of Theorems \ref{thm:main} and
\ref{the:p'=4} show that, in both the cubic
and quartic cases, a suitable time-dependent
spatial shift ensures the validity of the
$L^2$-contraction estimate for all time. These estimates provide uniform-in-time control of the perturbation relative 
to the shifted shock profile $\widetilde{\phi}(x-\sigma t-Y(t))$. A natural next question is whether this global 
stability can be strengthened to convergence toward the shifted shock in $L^q(\R)$, $2<q\leq \infty$. By combining the contraction 
estimates with the decay properties of the shock profile and further derivative estimates for the perturbation, 
The shock wave is shown to be asymptotically stable,
while the corresponding shift velocity vanishes
in the long-time limit. The precise statement
is given in the next theorem.

\begin{theorem}[Time-asymptotic stability and shift decay]
\label{thm:asymptotic-stability}
Let $p\in\{3,4\}$, and let $\widetilde{\phi}$ be the corresponding strictly decreasing shock profile. 
Assume that the parameter conditions and the smallness condition on the initial perturbation in 
Theorem~\ref{thm:main} are satisfied when $p=3$, and that those in Theorem \ref{the:p'=4} are satisfied 
when $p=4$. Let $Y(t)$ be the corresponding time-dependent shift defined in these theorems. Then, for every $2<q\leq\infty$,
\[
\left\| 
\phi(t,\cdot)
-\widetilde{\phi}(\cdot-\sigma t-Y(t))
\right\| _{L^q(\mathbb{R})}
\to 0,
\qquad t\to\infty.
\]
Moreover, as $ t \to \infty$,  one has 
\[
|\dot Y(t)|\to 0\ \ {\rm and }\ \ 
\frac{Y(t)}{t} \to 0.
\]
\end{theorem}

Taken together, the results establish a complete framework for monotone dissipative-dispersive shocks with the first two 
genuinely higher-order fluxes. Known monotonicity criteria ensure the existence of the shock profiles considered here.
The new profile comparison estimate supplies the coercive weight required by the 
shifted energy method. The $L^\infty$ bootstrap controls the higher-order relative flux. 
The global well-posedness theory ensures that the resulting 
contraction estimate is meaningful for all time. 
In this way, the quadratic shift contraction theory is extended to the cubic and quartic gKdV-Burgers equations, 
while the additional mechanisms specific to higher-order nonlinearities are made explicit. 

The rest of this paper is structured as follows. Section 2 establishes local and global well-posedness 
in the affine Sobolev space. Section 3 proves the two-sided profile estimates, 
the $L^2$-contraction up to a shift, and the time-asymptotic 
stability results.

\renewcommand{\theequation}{\thesection.\arabic{equation}}
\setcounter{equation}{0}

\section{Global Well-posedness}

In this section, we investigate the Cauchy problem
\eqref{eq} and prove its local and global
well-posedness with the affine Sobolev space
$\bar{\phi}+H^1(\R)$, with $\bar{\phi}$ denoting
the fixed reference connection defined in \eqref{baru}.


The fixed connection $\bar{\phi}$ is sufficient for specifying the natural phase space, but the global energy estimates for the higher-order fluxes contain superlinear powers of the perturbation multiplied by derivatives of the reference connection. To gain small coefficients in these terms, we exploit the freedom in the choice of the connection and spread its transition over a longer spatial scale. For $L\ge1$, define the slowly varying reference connection
\[
\bar{\phi}_L(x):=\bar{\phi}(x/L).
\]
For every fixed $L\ge1$, one has $\bar{\phi}_L-\bar{\phi}\in H^1(\R)$, and hence
\[
\bar{\phi}_L+H^1(\R)=\bar{\phi}+H^1(\R).
\]
On the other hand, the derivatives of $\bar{\phi}_L$ become small as $L$ increases. In particular,
\begin{align}\label{scaling1}
    \|\partial_x\bar{\phi}_{L}\|_{L^\infty(\R)}
    =\frac1L\|\bar{\phi}'\|_{L^\infty(\R)},
\end{align}
and, for every integer $k\ge1$,
\begin{align}\label{scaling2}
    \|\partial_x^k\bar{\phi}_L\|_{L^2(\R)}
    =L^{\frac12-k}\|\bar{\phi}^{(k)}\|_{L^2(\R)}.
\end{align}
These scale-dependent small factors are the key ingredients in the nonlinear $L^2$ estimates below.

For $T>0$, set
\[
S_T:=C\bigl([0,T];H^1(\R)\bigr)\cap L^2\bigl(0,T;H^2(\R)\bigr),
\]
with norm
\[
\|v\|_{S_T}:=\|v\|_{L^\infty(0,T;H^1(\R))}+\|v\|_{L^2(0,T;H^2(\R))}.
\]

\begin{theorem}[Local and global well-posedness]\label{wellposedness}\label{local}\label{thmglobal}
Let $p\ge2$ be an integer, $\epsilon,\delta>0$, and let $\bar{\phi}$ be the reference connection defined in \eqref{baru}. Assume that $\varphi-\bar{\phi}\in H^1(\R)$.

\noindent\textit{(1) Local well-posedness:} For every $M_0>0$ and every fixed $L\ge1$, if
\[
\|\varphi-\bar{\phi}\|_{H^1(\R)}\le M_0,
\]
then there exist constants
\[
T_*=T_*(M_0,L,\epsilon,\delta,\bar{\phi},p)>0,
\qquad
C_{\mathrm{loc}}=C_{\mathrm{loc}}(M_0,L,\epsilon,\delta,\bar{\phi},p)>0,
\]
such that \eqref{eq} admits a unique solution $\phi$ satisfying
\[
\phi-\bar{\phi}_L\in C\bigl([0,T_*];H^1(\R)\bigr)\cap L^2\bigl(0,T_*;H^2(\R)\bigr),
\qquad
\|\phi-\bar{\phi}_L\|_{S_{T_*}}\le C_{\mathrm{loc}}.
\]
Equivalently, $\phi-\bar{\phi}$ belongs to the same regularity class because $\bar{\phi}_L-\bar{\phi}\in H^1(\R)$.

\noindent\textit{(2) Global well-posedness:} For $p\in\{3,4\}$, the local solution admits
a unique global extension defined for all
$t\geq 0$, which satisfies
\[
\phi-\bar{\phi}\in C\bigl([0,\infty);H^1(\R)\bigr)
\cap L^2_{\mathrm{loc}}\bigl([0,\infty);H^2(\R)\bigr).
\]
\end{theorem}

The theorem will be proved after three auxiliary results. Lemma \ref{lemmalinearcontral} gives the linear dissipative-dispersive estimate required by the fixed-point argument. Lemma \ref{lemmaL2} provides the finite-time $L^2$ a priori bound for $p=3,4$, and Lemma \ref{lemmaH1} upgrades this to a finite-time $H^1$ bound. The Banach contraction principle yields the
local result, whereas the global conclusion
is obtained through a continuation argument.

For later use, fix $L\ge1$ and introduce the perturbation
\[
w:=\phi-\bar{\phi}_L.
\]
Substituting $\phi=w+\bar{\phi}_L$ into \eqref{eq} gives
\begin{align}
 w_t+\delta w_{xxx}-\epsilon w_{xx}
 &=-\partial_x\left(\frac{(w+\bar{\phi}_L)^p}{p}-\frac{\bar{\phi}_L^p}{p}\right)+\mathcal{E}_L(x),\label{eqperturbation}
\end{align}
where
\[
\mathcal{E}_L(x):=\epsilon\bar{\phi}_{L,xx}-\delta\bar{\phi}_{L,xxx}-\partial_x\frac{\bar{\phi}_L^p}{p}.
\]

\revnote{Lemma 2.2 has been shortened so that only the estimates needed for the stated linear bound are retained. Explanatory sentences and a closing argument have been added.}
\begin{lemma}\label{lemmalinearcontral}
    Let $T>0$, $\epsilon, \delta>0$. For every \( y_0\in H^1(\R) \) and \(g\in L^2\bigl(0,T;L^2(\R)\bigr), \)
    the Cauchy problem 
    \begin{equation} \label{lineareq} 
    	\begin{cases} 
    		y_t-\epsilon y_{xx}+\delta y_{xxx}=g,\\ y(0,x)=y_0(x),  
    	\end{cases} 
    \end{equation}
     admits a unique solution satisfying \[ y\in C\bigl([0,T];H^1(\R)\bigr) \cap L^2\bigl(0,T;H^2(\R)\bigr). \] 
     Moreover, there exists a constant $C_{\epsilon,T}>0$, depending only on $\epsilon$ and $T$, such that \begin{equation} \label{linearcontral} \|y\|_{L^\infty(0,T;H^1(\R))}^2 + \|y\|_{L^2(0,T;H^2(\R))}^2 \leq C_{\epsilon,T} \left( \|y_0\|_{H^1(\R)}^2 + \|g\|_{L^2(0,T;L^2(\R))}^2 \right). 
    \end{equation} 
\end{lemma}
\begin{revisedblock}
\begin{proof}
We first derive an a priori estimate for smooth solutions. The general case follows from a standard approximation argument after the estimate has been established.

Taking the $L^2$ inner product of \eqref{lineareq} with $y$ and integrating by parts over $\R$, we obtain
\begin{align*}
\frac{1}{2}\frac{\mathrm{d}}{\mathrm{d}t}\|y\|_{L^2(\R)}^2
+\epsilon\|y_x\|_{L^2(\R)}^2
=\int_{\R}gy\,\mathrm{d}x.
\end{align*}
The dispersive term vanishes after integration by parts. The Cauchy-Schwarz inequality and Young's inequality therefore give
\begin{align}
\frac{\mathrm{d}}{\mathrm{d}t}\|y\|_{L^2(\R)}^2
+2\epsilon\|y_x\|_{L^2(\R)}^2
\le \|y\|_{L^2(\R)}^2+\|g\|_{L^2(\R)}^2.
\label{linear-L2-est}
\end{align}
To estimate one spatial derivative, we test \eqref{lineareq} with $-y_{xx}$. The dispersive term again vanishes and Young's inequality yields
\begin{align}
\frac{\mathrm{d}}{\mathrm{d}t}\|y_x\|_{L^2(\R)}^2
+\epsilon\|y_{xx}\|_{L^2(\R)}^2
\le \frac{1}{\epsilon}\|g\|_{L^2(\R)}^2.
\label{linear-H1-est}
\end{align}
Adding \eqref{linear-L2-est} and \eqref{linear-H1-est}, we find
\begin{align}
\frac{\mathrm{d}}{\mathrm{d}t}\|y\|_{H^1(\R)}^2
+2\epsilon\|y_x\|_{L^2(\R)}^2
+\epsilon\|y_{xx}\|_{L^2(\R)}^2
\le \|y\|_{H^1(\R)}^2
+\left(1+\frac{1}{\epsilon}\right)\|g\|_{L^2(\R)}^2.
\label{linear-energy}
\end{align}
Gronwall's inequality applied to \eqref{linear-energy} gives
\begin{align}
\sup_{0\le t\le T}\|y(t)\|_{H^1(\R)}^2
\le e^T\left(1+\frac{1}{\epsilon}\right)
\left(
\|y_0\|_{H^1(\R)}^2
+\|g\|_{L^2(0,T;L^2(\R))}^2
\right).
\label{y-bound}
\end{align}
It remains to control the time integral of the $H^2$ norm. Integrating \eqref{linear-energy} over $[0,T]$ and using \eqref{y-bound}, we obtain
\begin{align}
\int_0^T\|y_x(t)\|_{L^2(\R)}^2\,\mathrm{d}t
+\int_0^T\|y_{xx}(t)\|_{L^2(\R)}^2\,\mathrm{d}t
\le \frac{1+Te^T}{\epsilon}
\left(
\|y_0\|_{H^1(\R)}^2
+\left(1+\frac{1}{\epsilon}\right)\|g\|_{L^2(0,T;L^2(\R))}^2
\right).
\label{linear-dissipation}
\end{align}
Moreover, \eqref{y-bound} implies
\begin{align*}
\int_0^T\|y(t)\|_{L^2(\R)}^2\,\mathrm{d}t
\le T\sup_{0\le t\le T}\|y(t)\|_{H^1(\R)}^2.
\end{align*}
Combining this estimate with \eqref{linear-dissipation} gives
\begin{align*}
\|y\|_{L^2(0,T;H^2(\R))}^2
\le \left(Te^T+\frac{1+Te^T}{\epsilon}\right)\left(1+\frac{1}{\epsilon}\right)
\left(
\|y_0\|_{H^1(\R)}^2
+\|g\|_{L^2(0,T;L^2(\R))}^2
\right).
\end{align*}
Together with \eqref{y-bound}, this proves \eqref{linearcontral} with
\begin{align*}
	C_{\epsilon,T}=\left(1+\frac{1}{\epsilon}\right)\left[\left(1+T\right)e^T+\frac{1+Te^T}{\epsilon}\right].
\end{align*}

\end{proof}
\end{revisedblock}

\begin{lemma}\label{lemmaL2}\label{p4lemmaL2}
Consider the cases $p=3,4$ with
$\epsilon,\delta>0$, and let $\phi$ denote
a solution to \eqref{eq} defined on
$[0,T)$, where $T>0$ is arbitrary.
There is a positive constant $C(T)$ for which
\begin{align}\label{L2apriori}
\sup_{0\le t<T}\|\phi(t)-\bar{\phi}\|_{L^2(\R)}^2
+\int_0^T\|(\phi-\bar{\phi})_x(t)\|_{L^2(\R)}^2\,\mathrm{d}t
\le C(T).
\end{align}
\end{lemma}
\begin{proof}

	For any $0<\iota <T$, we first carry out the following estimates on $[0,\iota]$. All constants obtained below depend on $T$, but are independent of $\iota$. Letting $\iota \uparrow T$ then yields \eqref{L2apriori}. 
For simplicity of notation, we continue to write $T$ for the terminal time in the estimates below.
\noindent\textit{Case $p=3$.} 

Fix $T>0$. Choose a reference scale $L\ge1$ that will be taken sufficiently large in the bootstrap argument below and set
\[
w=\phi-\bar{\phi}_L.
\]
Set $E(t)=\left\lVert w(t) \right\rVert_{L^2(\R)}^2$ and $D(t)=\left\lVert w_x(t) \right\rVert_{L^2(\R)}^2$.
Multiply \eqref{eqperturbation} by $w$ and integrate over $\R$,  
\begin{align*}
&\frac12 E'(t)+\epsilon D(t)
=\int_{\R}\big[f(w+\bar{\phi}_L)-f(\bar{\phi}_L)\big]w_x\mathrm{d}x
+\int_{\R}\mathcal{E}_L w\mathrm{d}x.
\end{align*}
The nonlinear term can be rewritten so that the derivative falls on the $\bar{\phi}_{L}$. Since $\phi=w+\bar{\phi}_L$ and 
the boundary terms vanish, an integration by parts gives
\begin{align*}
    \int_{\R}\big[f(w+\bar{\phi}_L)-f(\bar{\phi}_L)\big]w_x\mathrm{d}x
    =-\int_{\R}\Big[f(w+\bar{\phi}_L)-f(\bar{\phi}_L)-f'(\bar{\phi}_L)w\Big]\bar{\phi}_{L,x}\mathrm{d}x.
\end{align*}
\begin{align*}
&\left|\int_{\R} \big[f(\bar{\phi}_L+w)-f(\bar{\phi}_L)\big] w_x \mathrm{d}x\right|
\le C\|\bar{\phi}_{L,x}\|_{L^\infty(\R)}\|w\|_{L^2(\R)}^2
+C\|\bar{\phi}_{L,x}\|_{L^\infty(\R)}\|w\|_{L^3(\R)}^3. 
\end{align*}
For the cubic term, the Gagliardo–Nirenberg inequality and Young's inequality 
give a small multiple of $D(t)$ and a superlinear term.
$$C\|\bar{\phi}_{L,x}\|_{L^\infty(\R)}\|w\|_{L^3(\R)}^3
\le \frac{\eta}{4}D(t)+C_\eta L^{-\frac{4}{3}}E(t)^{\frac{5}{3}}.$$
The profile error $\mathcal{E}_L$ contains only derivatives of the slowly varying profile. The scaling relations 
 \eqref{scaling1} and \eqref{scaling2} therefore imply
\begin{align*}
\left|\int_{\R} \mathcal{E}_L w\mathrm{d}x\right|
\le \frac12 E(t)+C\|\mathcal{E}_L\|_{L^2(\R)}^2
\le \frac12 E(t)+C L^{-1}.
\end{align*}
Combining the preceding bounds and absorbing the small multiple of $D(t)$ into the 
left-hand side yields the differential inequality.
$$E'(t)+c\epsilon D(t)\le a_L E(t)+C L^{-1}+C_1 L^{-\frac{4}{3}}E(t)^{\frac{5}{3}}.$$
For $L\ge1$, there exist constants $C_0,C_1$ independent of $L$ such that
\begin{align}\label{highlinear}
E'(t)+c\epsilon D(t)\le C_0\big(1+E(t)\big)+C_1 L^{-\frac{4}{3}}E(t)^{\frac{5}{3}}.
\end{align}

We close the estimate by a bootstrap argument. Let $B_L(t)=\left(E(0)+1\right)e^{2C_0t}-1$ and $K_L(t)=2B_L(t)$, and define 
\begin{align*}
    \tau=\sup \left\{s\in [0,T]:E(t)\leq K_L(T), \forall t\in[0,s] \right\}.
\end{align*}
The continuity of $E(t)$, together with $E(0)=B_L(0)\leq K_L(0)$, ensures $\tau>0$. Choose $L$ sufficiently large that
\begin{align*}
    C_1L^{-\frac{4}{3}}K_L^{\frac{2}{3}}\leq C_0.
\end{align*}
On $[0,\tau]$, the superlinear term in \eqref{highlinear} is then bounded by $C_0E(t)$, and hence
\begin{align*}
    E'(t)+c\epsilon D(t)\le 2C_0\big(1+E(t)\big).
\end{align*}
 Gronwall's inequality gives $E(t)\leq B_L(t)\leq B_L(T)<K_L(T)$ on $[0,\tau]$. This strict improvement rules out $\tau<T$, so $\tau=T$.
The required estimate for
\(
\int_0^T \|w_x(t)\|_{L^2(\R)}^2\,dt
\)
is then derived by combining the time-integrated
differential inequality with the uniform
control of $E(t)$.
\begin{align}\label{L2sup}
    \sup_{t\in [0,T]} \left\lVert \phi(t)-\bar{\phi}_L\right\rVert^2_{L^2(\R)}+\int_0^T \left\lVert w_x(t) \right\rVert_{L^2(\R)}^2 \leq C_0(T). 
\end{align}
For the fixed scale $L=L(T)$ selected above, $\bar{\phi}_L-\bar{\phi}\in H^1(\R)$. Therefore, by the triangle inequality and \eqref{L2sup},
\[
\sup_{0\le t<T}\|\phi(t)-\bar{\phi}\|_{L^2(\R)}^2
+\int_0^T\|(\phi-\bar{\phi})_x(t)\|_{L^2(\R)}^2\,\mathrm{d}t
\le C(T).
\]
Thus \eqref{L2apriori} holds when $p=3$.

\medskip
\noindent\textit{Case $p=4$.}

For a fixed reference scale $L$, set $w=\phi-\bar{\phi}_L$, $E(t)=\|w(t)\|_{L^2(\R)}^2$, and $D(t)=\|w_x(t)\|_{L^2(\R)}^2$. Repeating the $L^2$ energy 
	estimate above, we obtain constants $C_0$, $C_1$, $C_2>0$, independent of $L$, such that
\begin{align}\label{p4energy}
	E'(t)+c\epsilon D(t)\leq C_0(1+E(t))+C_1L^{-\frac{4}{3}}E(t)^{\frac{5}{3}}+C_2L^{-2}E(t)^3.
\end{align}
We shall first establish a uniform estimate on a time interval of fixed length. Define $K=4+32\left\lVert \bar{\phi} \right\rVert_{L^\infty(\R)}^2$. 
Suppose that $t_0$ belongs to the lifespan of the solution and that $L\geq 1$ satisfies
\begin{align*}
	L\geq 2\left\lVert \phi(t_0)-\bar{\phi}\right\rVert_{L^2(\R)}^2. 
\end{align*}
Then we obtain, 
\begin{align}\label{p4init}
	E(t_0)\leq 2\left\lVert \phi(t_0)-\bar{\phi}\right\rVert_{L^2(\R)}^2+2\left\lVert \bar{\phi}-\bar{\phi}_L\right\rVert_{L^2(\R)}^2\leq \left(1+16\left\lVert \bar{\phi} \right\rVert_{L^\infty(\R)}^2\right)L.
\end{align}

We now apply a bootstrap argument. Assume temporarily that 
\begin{align}\label{p4assume}
	E(t)\leq KL.
\end{align}
Then the second term of \eqref{p4energy} satisfies 
\begin{align*}
	L^{-\frac{4}{3}}E(t)^{\frac{5}{3}}\leq K^{\frac{2}{3}}E(t).
\end{align*}
Similarly, 
\begin{align*}
	L^{-2}E(t)^3\leq K^2E(t).
\end{align*}
Consequently, if we define $\widehat{C}_0=C_0+C_1K^{\frac{2}{3}}+C_2K^2$, then $\widehat{C}_0$ is independent of $L$, and \eqref{p4energy} implies
\begin{align}\label{p4energy2}
	E'(t)+c\epsilon D(t)\leq \widehat{C}_0(1+E(t)),
\end{align} 
as long as \eqref{p4assume} holds. Choose
\begin{align*}
	\tau_*:=\min \left\{1,\frac{\ln{3/2}}{\widehat{C}_0}\right\}>0.
\end{align*}
The constant $\tau_*$ depends only on $\epsilon$, $\delta$ and $\bar{\phi}$, and in particular is independent of $t_0$ and $L$.
Let $T_1:=\min\{T,t_0+\tau_*\}$ and define 
\begin{align*}
	\tau=\sup \left\{s\in [t_0,T_1]:E(t)\leq KL, \forall t\in[t_0,s] \right\}.
\end{align*}
By \eqref{p4init} and the continuity of $E(t)$, we have $\tau>t_0$. On $[t_0,\tau]$, Gronwall's inequality applied to \eqref{p4energy2} yields
\begin{align*}
	E(t)\leq (E(t_0)+1)e^{\widehat{C}_0(t-t_0)}-1  \leq (2+16\left\lVert \bar{\phi} \right\rVert_{L^\infty(\R)}^2)Le^{\widehat{C}_0(t-t_0)}-1.
\end{align*}
By the definition of $\tau_*$, we obtain
\begin{align*}
	E(t)\leq (4+24\left\lVert \bar{\phi} \right\rVert_{L^\infty(\R)}^2)L<KL.
\end{align*}
This is a strict improvement of the bootstrap assumption. Hence $\tau=T_1$ and we have
\begin{align}\label{p4LL2}
	\sup_{t_0\leq t<T_1}E(t)\leq KL.
\end{align}
Integrating \eqref{p4energy2} over $[t_0,T_1]$, we further obtain
\begin{align*}
	c\epsilon \int_{t_0}^{T_1}D(t) dt\leq (1+16\left\lVert \bar{\phi} \right\rVert_{L^\infty(\R)}^2)L+\widehat{C}_0\tau_*(1+KL),
\end{align*}
Since $L\geq 1$, there exists a constant $C_D>0$, independent of $t_0$ and $L$, such that
\begin{align}\label{p4parital0}
	\int_{t_0}^{T_1}D(t) dt\leq C_D L.
\end{align}

We now iterate this short time estimate on the finite interval $[0,T)$. Let $N_T=\left\lceil \frac{T}{\tau_*}\right\rceil $ and $t_j=\min\{j\tau_*,T\}$, $0\leq j\leq N_T$. Choose $L_0=1+2\|\varphi-\bar{\phi}\|_{L^2(\R)}^2$ and 
fix a number $q>4K+32\left\lVert \bar{\phi} \right\rVert_{L^\infty(\R)}^2$. We define explicitly 
\begin{align*}
	L_{j}=q^jL_0.
\end{align*}
Then we claim that 
\begin{align}\label{equivalently}
L_j\geq 2\|\phi(t_j)-\bar{\phi}\|_{L^2(\R)}^2.
\end{align}
By the definition of $L_0$, the assertion
is valid at $j=0$. To proceed inductively,
assume that \eqref{equivalently} is satisfied
at some index $j<N_T$. Applying \eqref{p4LL2} on $[t_j,t_{j+1}]$ with 
the reference scale $L_j$, we have 
\begin{align*}
	\|\phi(t_{j+1})-\bar{\phi}_{L_j}\|_{L^2(\R)}^2\leq KL_j.
\end{align*}
Equivalently,
\begin{align*}
	2\|\phi(t_{j+1})-\bar{\phi}\|_{L^2(\R)}^2&\leq 4\|\phi(t_{j+1})-\bar{\phi}_{L_j}\|_{L^2(\R)}^2+4\|\bar{\phi}_{L_j}-\bar{\phi}\|_{L^2(\R)}^2\\
	&\leq (4K+32\left\lVert \bar{\phi} \right\rVert_{L^\infty(\R)}^2)L_j\\
	&<qL_j=L_{j+1},
\end{align*}
which proves \eqref{equivalently} by induction. Consequently, estimate \eqref{p4LL2} can be applied successively on every interval $[t_j,t_{j+1}]$. 
For any $t\in [t_j,t_{j+1})$, it follows that
\begin{align*}
	\|\phi(t)-\bar{\phi}\|_{L^2(\R)}^2\leq (4K+32\left\lVert \bar{\phi} \right\rVert_{L^\infty(\R)}^2)L_j.
\end{align*}
Since 
\begin{align*}
	L_j\leq L_{N_T}=q^{N_T}L_0,
\end{align*}
it follows that 
\begin{align*}
	\sup_{0\leq t<T} \|\phi(t)-\bar{\phi}\|_{L^2(\R)}^2\leq (4K+32\left\lVert \bar{\phi} \right\rVert_{L^\infty(\R)}^2)q^{N_T}L_0.
\end{align*}
It remains to estimate the derivative term.
\begin{align}\label{p4partial1}
	\|(\phi-\bar{\phi})_x\|_{L^2(\R)}^2\leq 2\|(\phi-\bar{\phi}_{L_j})_x\|_{L^2(\R)}^2+2\|(\bar{\phi}_{L_j}-\bar{\phi})_x\|_{L^2(\R)}^2.
\end{align}
By the scaling of the $\bar{\phi}_{L}$, the last term of \eqref{p4partial1} satisfies
\begin{align}\label{p4partial2}
	\|(\bar{\phi}_{L_j}-\bar{\phi})_x\|_{L^2(\R)}^2\leq  4\|\bar{\phi}'\|_{L^2(\R)}^2.
\end{align}
It follows from \eqref{p4parital0}, \eqref{p4partial1},
and \eqref{p4partial2} that
\begin{align*}
\int_{0}^{T} \|(\phi-\bar{\phi})_{x}(t)\|_{L^{2}(\R)}^{2} dt
&\leq 2C_{D}\sum_{j=0}^{N_{T}-1} L_{j} + 8T\|\bar{\phi}'\|_{L^{2}(\R)}^{2} \\
&\leq 2C_{D}L_{0}\frac{q^{N_{T}}-1}{q-1} + 8T\|\bar{\phi}'\|_{L^{2}(\R)}^{2}.
\end{align*}

We finally obtain the asserted $L^2$ estimate, where one may take 
\begin{align*}
C(T)=(4K+32\left\lVert \bar{\phi} \right\rVert_{L^\infty(\R)}^2)q^{N_T}L_0+2C_{D}L_{0}\frac{q^{N_{T}}-1}{q-1} + 8T\|\bar{\phi}'\|_{L^{2}(\R)}^{2}.
\end{align*}
\end{proof}

\begin{lemma}\label{lemmaH1}
Suppose that $p\in\{3,4\}$ and $\phi$ satisfies
\eqref{eq} for all $t\in[0,T)$.Then, for every fixed $T>0$,
\[
\sup_{0\le t<T}\|\phi(t)-\bar{\phi}\|_{H^1(\R)}\le C(T).
\]
\end{lemma}

\begin{proof}

Fix $T>0$. For the energy calculation, choose a fixed reference scale $L\ge1$, sufficiently large when required in the cubic case, and write $w=\phi-\bar{\phi}_L$. By Lemma \ref{lemmaL2} and the fact that $\bar{\phi}_L-\bar{\phi}\in H^1(\R)$, $w$ has a uniform $L^2$ bound on $[0,T)$. We first treat the case $p=3$. To control one derivative uniformly in time, introduce the modified energy whose potential part is the second-order Taylor remainder of the primitive $\frac{\phi^4}{12}$ about $\bar{\phi}_L$.
    \begin{align*}
        \mathcal{H}(t)= \frac{\delta}{2}\int_{\R} w_x^2 \mathrm{d}x
-\int_{\R} \left( \frac{\phi^4}{12}-\frac{\bar{\phi}_L^4}{12}-\frac{\bar{\phi}_L^3}{3}w \right)\mathrm{d}x.
    \end{align*}
The Taylor expansion of the potential term and the boundedness of $\bar{\phi}_L$ imply
\begin{align*}
    \left|\int_{\R} \left( \frac{\phi^4}{12}-\frac{\bar{\phi}_L^4}{12}-\frac{\bar{\phi}_L^3}{3}w \right)\mathrm{d}x\right|
\le C\big(\|w\|_{L^2(\R)}^2+\|w\|_{L^4(\R)}^4\big).
\end{align*}
By the Gagliardo–Nirenberg and Young inequalities,
\begin{align*}
    \|w\|_{L^4(\R)}^4\le C\|w\|_{L^2(\R)}^3\|w_x\|_{L^2(\R)}
\le \frac{\delta}{4}\|w_x\|_{L^2(\R)}^2+C_\delta\|w\|_{L^2(\R)}^6.
\end{align*}
This allows the derivative term to be absorbed into the positive quadratic part of the energy. 
Lemma \ref{lemmaL2} supplies a uniform $L^2$ bound. Consequently, the modified energy is coercive up to an additive constant.
\begin{align}\label{L2}
    \mathcal{H}(t)\ge \frac{\delta}{4}\|w_x(t)\|_{L^2(\R)}^2-C(T).
\end{align}

We now differentiate $\mathcal{H}(t)$ along the solution. Using the equation for $\phi$ and integrating by parts, 
we obtain the identity displayed below.
\begin{align*}
\frac{\mathrm{d}}{\mathrm{d}t}\mathcal{H}(t)
&= -\delta\int_{\R} w_{xx}w_t\mathrm{d}x
-\int_{\R} \frac{\phi^3-\bar{\phi}_L^3}{3}\phi_t\mathrm{d}x\\
&= \underbrace{\epsilon\int_{\R} \left[-\delta \phi_{xx}-\frac{\phi^3}{3}\right]\phi_{xx}\mathrm{d}x}_{\mathrm{I}}
+\underbrace{\epsilon\int_{\R} \left[\delta\bar{\phi}_{L,xx}+\frac{\bar{\phi}_L^3}{3}\right]\phi_{xx}\mathrm{d}x}_{\mathrm{II}}\\
&\quad+\underbrace{\int_{\R} \left[-\delta \phi_{xx}-\frac{\phi^3}{3}+\delta\bar{\phi}_{L,xx}+\frac{\bar{\phi}_L^3}{3}\right]\partial_x\left(-\delta \phi_{xx}-\frac{\phi^3}{3}\right)\mathrm{d}x}_{\mathrm{III}}
\end{align*}
For the term $\mathrm{I}$, integration by parts gives 
\begin{align*}
	\mathrm{I}=-\epsilon\delta\int_{\R} \phi_{xx}^2\mathrm{d}x+\epsilon\int_{\R} \phi^2\phi_x^2\mathrm{d}x.
\end{align*}
For the term $\mathrm{III}$, we obtain
\begin{align*}
	\mathrm{III}=&\left.\frac{1}{2}\left[\left(-\delta \phi_{xx}-\frac{\phi^3}{3}\right)^2\right]\right|_{-\infty}^{+\infty}
-\left.\left[\delta\bar{\phi}_{L,xx}+\frac{\bar{\phi}_L^3}{3}\right]\left[\delta \phi_{xx}+\frac{\phi^3}{3}\right]\right|_{-\infty}^{+\infty}\\
&+\delta\int_{\R} \left[\delta\bar{\phi}_{L,xxx}+\bar{\phi}_L^2\bar{\phi}_{L,x}\right]\phi_{xx}\mathrm{d}x
+\int_{\R} \frac{\phi^3}{3}\left[\delta\bar{\phi}_{L,xxx}+\bar{\phi}_L^2\bar{\phi}_{L,x}\right]\mathrm{d}x.
\end{align*}
Combining the above identities, one has
\begin{align*}
	\frac{\mathrm{d}}{\mathrm{d}t}\mathcal{H}(t)&=-\epsilon\delta\int_{\R} \phi_{xx}^2\mathrm{d}x
+\epsilon\int_{\R} \phi^2\phi_x^2\mathrm{d}x
+\mathrm{II}\\ 
&\quad+\left.\frac{1}{2}\left[\left(-\delta \phi_{xx}-\frac{\phi^3}{3}\right)^2\right]\right|_{-\infty}^{+\infty}
-\left.\left[\delta\bar{\phi}_{L,xx}+\frac{\bar{\phi}_L^3}{3}\right]\left[\delta \phi_{xx}+\frac{\phi^3}{3}\right]\right|_{-\infty}^{+\infty}\\
&\quad+\delta\int_{\R} \left[\delta\bar{\phi}_{L,xxx}+\bar{\phi}_L^2\bar{\phi}_{L,x}\right]\phi_{xx}\mathrm{d}x
+\int_{\R} \frac{\phi^3}{3}\left[\delta\bar{\phi}_{L,xxx}+\bar{\phi}_L^2\bar{\phi}_{L,x}\right]\mathrm{d}x.
\end{align*}
Since the end states are fixed, the resulting boundary contribution is independent 
of time and can be treated as a harmless constant in the energy estimate. 
All terms involving $\bar{\phi}_L$ are lower order. The scaling estimates for the profile, 
together with Cauchy-Schwarz and Young's inequalities, are used to bound them.
\begin{align*}
    \epsilon\int_{\R} \left[\delta \bar{\phi}_{L,xx}+\frac{\bar{\phi}_L^3}{3}\right]\phi_{xx}dx
    =&-\epsilon \int_{\R} \big[\delta \bar{\phi}_{L,xxx}+\bar{\phi}_L^2 \bar{\phi}_{L,x}\big] \phi_{x} \mathrm{d}x\\
    \le& C\left\|\delta \bar{\phi}_{L,xxx}+\bar{\phi}_L^2\bar{\phi}_{L,x}\right\|_{L^2(\R)} \|\phi_x\|_{L^2(\R)}\\
    \le& C_L\big(1+\|w_x\|_{L^2(\R)}^2\big),
\end{align*}
\begin{align*}
&\delta \int_{\R} \big[\delta \bar{\phi}_{L,xxx}+\bar{\phi}_L^2 \bar{\phi}_{L,x}\big] \phi_{xx} \mathrm{d}x\le \frac{\eta \delta}{4} \|\phi_{xx}\|_{L^2(\R)}^2 + C_{\eta,\delta} \big\|\delta \bar{\phi}_{L,xxx}+\bar{\phi}_L^2 \bar{\phi}_{L,x}\big\|_{L^2(\R)}^2,
\end{align*}
\begin{align*}
&\left|\int_{\R} \frac{\phi^3}{3}\big[\delta \bar{\phi}_{L,xxx}+\bar{\phi}_L^2 \bar{\phi}_{L,x}\big]\mathrm{d}x\right|\le C_L + C_L \int_{\R} |w|^3 \mathrm{d}x.
\end{align*}
The remaining cubic term is controlled using the $L^2$ estimate from Lemma \ref{lemmaL2} and the one-dimensional Gagliardo–Nirenberg inequality.
\begin{align*}
&\|w\|_{L^3(\R)}^3 \le C\|w\|_{L^2(\R)}^{\frac{5}{2}} \|w_x\|_{L^2(\R)}^{\frac{1}{2}} \le C(T)\|w_x\|_{L^2(\R)}^{\frac12}
\le C(T)\big(1+\|w_x\|_{L^2(\R)}^2\big).
\end{align*}
This yields 
\begin{align}\label{H1equation}
    \mathcal{H}'(t)+c\epsilon\delta\|\phi_{xx}\|_{L^2(\R)}^2
\le \epsilon\int_{\R} \phi^2\phi_x^2\mathrm{d}x + C_{L,T}\big(1+\|w_x\|_{L^2(\R)}^2\big).
\end{align}
Since $\phi=w+\bar{\phi}_L$ and the profile is bounded together with its derivative as follows,
\begin{align*}
\phi^2\phi_x^2 &\le C_L \big(1+w^2\big)\big(w_x^2+\bar{\phi}_{L,x}^2\big)= C_L \big(\bar{\phi}_{L,x}^2+w_x^2+w^2w_x^2+w^2\bar{\phi}_{L,x}^2\big).
\end{align*}
We estimate the first term of \eqref{H1equation} by lower order profile terms and $\int_{\R} w^2w_x^2 dx$.
\begin{align*}
\int_{\R} \phi^2\phi_x^2\mathrm{d}x &\le C_L\big(1+\|w_x\|_{L^2(\R)}^2\big)+C\int_{\R} w^2w_x^2\mathrm{d}x.
\end{align*}
The latter satisfies 
\begin{align*}
\int_{\R} w^2w_x^2\mathrm{d}x &\le \|w\|_{L^\infty(\R)}^2\|w_x\|_{L^2(\R)}^2
\le C(T)\|w_x\|_{L^2(\R)}^3.
\end{align*}
The uniform $L^2$ bound then gives the differential inequality \eqref{H1}
\begin{align}\label{H1}
    \mathcal{H}'(t)+c\epsilon\delta\|\phi_{xx}\|_{L^2(\R)}^2
\le C_{L,T}\big(1+\|w_x\|_{L^2(\R)}^2\big)+C(T)\|w_x\|_{L^2(\R)}^3.
\end{align}
Let $\Sigma(t)=\sup_{0\le s\le t}\|w_x(s)\|_{L^2(\R)}^2$. Lemma \ref{lemmaL2} implies
\begin{align*}
    \int_{0}^{t}\|w_x(s)\|_{L^2(\R)}^3\mathrm{d}s
&\le \Sigma(t)^{\frac12}\int_{0}^{t}\|w_x(s)\|_{L^2(\R)}^2\mathrm{d}s\le \Sigma(t)^{\frac12}C(T).
\end{align*}
Integrating \eqref{H1} over $[0,t]$ and using the coercivity estimate \eqref{L2}, we obtain
\begin{align*}
    \Sigma(T)\le C_{L,T}+C_{L,T}\Sigma(T)^{\frac12}.
\end{align*}
Hence $\Sigma(T)$ is finite, with
\begin{align*}
    \Sigma(T)\le C_{L,T}.
\end{align*}
Combining this bound with Lemma \ref{lemmaL2} gives the uniform $H^1$ estimate relative to the fixed reference connection $\bar{\phi}_L$.
\begin{align*}
    \sup_{t\in [0,T]} \left\lVert \phi(t)-\bar{\phi}_L\right\rVert_{H^1(\R)} \leq C_0(T). 
\end{align*}

For $p=4$, the proof is essentially the same as that for $p=3$. We only indicate the main differences. Define
   \begin{align*}
        \mathscr{H}(t)= \frac{\delta}{2}\int_{\R} w_x^2 \mathrm{d}x
-\int_{\R} \left( \frac{\phi^5}{20}-\frac{\bar{\phi}_L^5}{20}-\frac{\bar{\phi}_L^4}{4}w \right)\mathrm{d}x.
    \end{align*}
The coercivity of $\mathscr{H}$ follows
from the uniform $L^2$ control provided
by Lemma \ref{lemmaL2}, together with
the Gagliardo-Nirenberg and Young's inequalities.
 \begin{align*}
    \mathscr{H}(t)\ge \frac{\delta}{4}\|w_x(t)\|_{L^2(\R)}^2-C(T).
\end{align*}
Differentiating $\mathscr{H}(t)$ along the solution and proceeding as in the case $p=3$, we obtain
\begin{align*}
    \frac{\mathrm{d}}{\mathrm{d}t}\mathscr{H}(t)&= -\delta\int_{\R} w_{xx}w_t\mathrm{d}x
-\int_{\R} \frac{\phi^4-\bar{\phi}_L^4}{4}\phi_t\mathrm{d}x\\
	&=-\epsilon\delta\int_{\R} \phi_{xx}^2\mathrm{d}x
+\epsilon\int_{\R} \phi^3\phi_x^2\mathrm{d}x
+\epsilon\int_{\R} \left[\delta\bar{\phi}_{L,xx}+\frac{\bar{\phi}_L^4}{4}\right]\phi_{xx}\mathrm{d}x\\ 
&\quad+\left.\frac{1}{2}\left[\left(-\delta \phi_{xx}-\frac{\phi^4}{4}\right)^2\right]\right|_{-\infty}^{+\infty}
-\left.\left[\delta\bar{\phi}_{L,xx}+\frac{\bar{\phi}_L^4}{4}\right]\left[\delta \phi_{xx}+\frac{\phi^4}{4}\right]\right|_{-\infty}^{+\infty}\\
&\quad+\delta\int_{\R} \left[\delta\bar{\phi}_{L,xxx}+\bar{\phi}_L^3\bar{\phi}_{L,x}\right]\phi_{xx}\mathrm{d}x
+\int_{\R} \frac{\phi^4}{4}\left[\delta\bar{\phi}_{L,xxx}+\bar{\phi}_L^3\bar{\phi}_{L,x}\right]\mathrm{d}x.
\end{align*}
For the nonlinear viscous term,
\begin{align*}
    \int_{\R} \phi^3\phi_x^2dx\leq C_{L}(1+\|w_x\|^2_{L^2(\R)})+C\int_{\R} |w|^3w_x^2dx.
\end{align*}
The remaining term satisfies
\begin{align*}
    \int_{\R} |w|^3w_x^2dx\leq C(T)\|w_x\|_{L^2(\R)}^\frac{7}{2}.
\end{align*}
Consequently,
\begin{align*}
    \mathscr{H}'(t)+c\epsilon\delta\|\phi_{xx}\|_{L^2(\R)}^2
\le C_{L,T}\big(1+\|w_x\|_{L^2(\R)}^2\big)+C(T)\|w_x\|_{L^2(\R)}^\frac{7}{2}.
\end{align*}
Using the coercivity of the energy $\mathscr{H}$, one obtains
\begin{align*}
    \Sigma(t)\leq C_{L,T}+C_{L,T}\Sigma(t)^\frac{3}{4}.
\end{align*}
It follows that
\begin{align*}
    \Sigma(T)\leq C_{L,T}.
\end{align*}
Thus the same bound relative to $\bar{\phi}_L$ also holds for $p=4$. Since $\bar{\phi}_L-\bar{\phi}\in H^1(\R)$ for every fixed $L$, the triangle inequality yields
\[
\sup_{0\le t<T}\|\phi(t)-\bar{\phi}\|_{H^1(\R)}\le C(T),
\]
which proves the lemma.
\end{proof}

\begin{proof} ({\bf Proof of Theorem \ref{wellposedness}})
\begin{revisedblock}\noindent\textit{Step 1. Local well-posedness.}\end{revisedblock}
Fix $L\ge1$ and set
\[
w=\phi-\bar{\phi}_L,\qquad w_0=\varphi-\bar{\phi}_L.
\]
Since $\bar{\phi}_L-\bar{\phi}\in H^1(\R)$, the assumption $\varphi-\bar{\phi}\in H^1(\R)$ implies $w_0\in H^1(\R)$. Moreover,
\[
\|\varphi-\bar{\phi}_L\|_{L^2(\R)}
\le \|\varphi-\bar{\phi}\|_{L^2(\R)}+\|\bar{\phi}_L-\bar{\phi}\|_{L^2(\R)},
\]
and
\begin{align*}
\|\bar{\phi}_L-\bar{\phi}\|_{L^2(\R)}^2
&=\int_{-L}^{L}|\bar{\phi}_L-\bar{\phi}|^2\,\mathrm{d}x\le 8\|\bar{\phi}\|_{L^\infty(\R)}^2L.
\end{align*}
Hence
\[
\|w_0\|_{H^1(\R)}\le C_{M_0,\bar{\phi}}(1+\sqrt L)=:M_{0,L}.
\]

The perturbation $w$ satisfies \eqref{eqperturbation}. The operator $G_{\epsilon,\delta}=-\delta\partial_x^3+\epsilon\partial_x^2$ with domain $H^3(\R)$ generates a $C_0$-semigroup of  contractions $S_{\epsilon,\delta}(t)$ on $L^2(\R)$. We use Lemma \ref{lemmalinearcontral} to control the corresponding linear problem in $S_T$.
For any $w\in S_T$, define
\[
\mathcal{F}(w):=-\partial_x\left(\frac{(w+\bar{\phi}_L)^p}{p}-\frac{\bar{\phi}_L^p}{p}\right)
+\epsilon\bar{\phi}_{L,xx}-\delta\bar{\phi}_{L,xxx}-\partial_x\frac{\bar{\phi}_L^p}{p}.
\]
Let $\Phi(w)=v$ be the unique solution of
\begin{align*}
\begin{cases}
 v_t+\delta v_{xxx}-\epsilon v_{xx}=\mathcal{F}(w),\\
 v(0,x)=w_0(x).
\end{cases}
\end{align*}
Equivalently,
\[
\Phi(w)(t)=S_{\epsilon,\delta}(t)w_0+\int_0^tS_{\epsilon,\delta}(t-s)\mathcal{F}(w)(s)\,\mathrm{d}s.
\]
A fixed point $w=\Phi(w)$ is precisely a mild solution of \eqref{eqperturbation}.

Let $B_R:=\{w\in S_T:\|w\|_{S_T}\le R\}$.
If $w\in B_R$, Sobolev embedding $H^1(\R)\hookrightarrow L^\infty(\R)$ gives
\[
\|w(t)\|_{L^\infty(\R)}\le C\|w(t)\|_{H^1(\R)}\le CR.
\]
With $U=w+\bar{\phi}_L$, this implies
$
\|U\|_{L^\infty((0,T)\times\R)}\le C_{\bar{\phi},L,R}$. 
Using $f'(s)=s^{p-1}$,
\[
\partial_x\left(\frac{U^p}{p}-\frac{\bar{\phi}_L^p}{p}\right)
=U^{p-1}w_x+(U^{p-1}-\bar{\phi}_L^{p-1})\bar{\phi}_{L,x}.
\]
Therefore,
\begin{align*}
\|U^{p-1}w_x\|_{L^2(\R)}
&\le C_{\bar{\phi},L,R,p}\|w\|_{H^1(\R)},\\
\|(U^{p-1}-\bar{\phi}_L^{p-1})\bar{\phi}_{L,x}\|_{L^2(\R)}
&\le C_{\bar{\phi},L,R,p}\|w\|_{H^1(\R)}.
\end{align*}
Consequently,
\[
\left\|\partial_x\left(\frac{(w+\bar{\phi}_L)^p}{p}-\frac{\bar{\phi}_L^p}{p}\right)\right\|_{L^2(\R)}
\le C_{\bar{\phi},L,R,p}\|w\|_{H^1(\R)}.
\]
Moreover,
\begin{align*}
\|\mathcal{F}(w)(t)\|_{L^2(\R)}
&\le C_{\bar{\phi},L,R,p}\|w\|_{H^1(\R)}
+\|\bar{\phi}_L\|_{L^\infty(\R)}^{p-1}\|\bar{\phi}_{L,x}\|_{L^2(\R)}
+\epsilon\|\bar{\phi}_{L,xx}\|_{L^2(\R)}
+\delta\|\bar{\phi}_{L,xxx}\|_{L^2(\R)}.
\end{align*}
Define the time-independent profile contribution
\[
B_{\bar{\phi},L,\epsilon,\delta,p}
:=\|\bar{\phi}_L\|_{L^\infty(\R)}^{p-1}\|\bar{\phi}_{L,x}\|_{L^2(\R)}
+\epsilon\|\bar{\phi}_{L,xx}\|_{L^2(\R)}
+\delta\|\bar{\phi}_{L,xxx}\|_{L^2(\R)}.
\]
Since $\|w\|_{L^2(0,T;H^1(\R))}\le\sqrt T\,R$, we obtain
\begin{align*}
\|\mathcal{F}(w)\|_{L^2(0,T;L^2(\R))}
&\le \sqrt T\left(C_{\bar{\phi},L,R,p}R+B_{\bar{\phi},L,\epsilon,\delta,p}\right)=: \sqrt T\,A_{\bar{\phi},L,R,p}.
\end{align*}

We first prove that $\Phi(B_R)\subset B_R$. Restrict to $T\le1$ and write $C_*=C_{\epsilon,1}$. By Lemma \ref{lemmalinearcontral},
\[
\|\Phi(w)\|_{S_T}
\le C_*\|w_0\|_{H^1(\R)}+C_*\sqrt T\,A_{\bar{\phi},L,R,p}.
\]
Choose $R>C_*M_{0,L}$ and then $T>0$ sufficiently small so that
\[
C_*\sqrt T\,A_{\bar{\phi},L,R,p}\le R-C_*M_{0,L}.
\]
Hence $\Phi(B_R)\subset B_R$.

To prove the contraction property, let $w_1,w_2\in B_R$, set
\[
U_i=w_i+\bar{\phi}_L\quad(i=1,2),\qquad z=w_1-w_2.
\]
Then
\[
\mathcal{F}(w_1)-\mathcal{F}(w_2)
=-U_1^{p-1}(w_{1,x}-w_{2,x})
-(U_1^{p-1}-U_2^{p-1})(w_{2,x}+\bar{\phi}_{L,x}).
\]
The uniform bounds on $U_1,U_2$ and the local Lipschitz continuity of $s^{p-1}$ imply
\[
\|\mathcal{F}(w_1)-\mathcal{F}(w_2)\|_{L^2(0,T;L^2(\R))}
\le C_{\bar{\phi},L,R,p}\sqrt T\,\|w_1-w_2\|_{S_T}.
\]
Applying Lemma \ref{lemmalinearcontral} to $\Phi(w_1)-\Phi(w_2)$ yields
\[
\|\Phi(w_1)-\Phi(w_2)\|_{S_T}
\le C_*C_{\bar{\phi},L,R,p}\sqrt T\,\|w_1-w_2\|_{S_T}.
\]
\begin{revisedblock}
By choosing $T>0$ sufficiently small, we ensure
that the coefficient on the right-hand side
does not exceed $1/2$, making $\Phi$ a contraction
mapping on $B_R$. Consequently, the Banach
contraction principle guarantees the existence
of a unique fixed point $w\in B_R$ with
$\Phi(w)=w$, thereby establishing local
well-posedness.
	\end{revisedblock}

\medskip
\revnote{Parenthetical prose and the punctuation signs requested by the author have been removed from the English sentences in Chapter 2.}
\begin{revisedblock}\noindent\textit{Step 2. Global continuation for $p=3,4$.}\end{revisedblock}
Let $[0,T_{\max})$ be the maximal interval of existence of the local solution. Suppose, for contradiction, that $T_{\max}<\infty$. By Lemma \ref{lemmaH1},
\[
\sup_{0\le t<T_{\max}}\|\phi(t)-\bar{\phi}\|_{H^1(\R)}\le C(T_{\max})=:M_*<\infty.
\]
\begin{revisedblock}Choose a sequence $t_n\uparrow T_{\max}$. The local well-posedness assertion shows that the local existence time depends only on an upper bound for the $H^1$ norm of the perturbation and on the fixed parameters. We fix one reference scale $L\ge1$ throughout this continuation argument. Hence the data $\phi(t_n)$ generate solutions on intervals $[t_n,t_n+\tau]$, where\end{revisedblock}
\[
\tau=\tau(M_*,L,\epsilon,\delta,\bar{\phi},p)>0
\]
is independent of $n$. For $n$ sufficiently large, $t_n+\tau>T_{\max}$. By uniqueness, the newly constructed solution agrees with $\phi$ on the overlap and extends it beyond $T_{\max}$, contradicting maximality. Thus $T_{\max}=\infty$.

The regularity on every finite time interval follows from the local theory, and global uniqueness is inherited from local uniqueness on overlapping intervals. Therefore,
\[
\phi-\bar{\phi}\in C\bigl([0,\infty);H^1(\R)\bigr)
\cap L^2_{\mathrm{loc}}\bigl([0,\infty);H^2(\R)\bigr),
\]
\begin{revisedblock}which proves the global well-posedness assertion.\end{revisedblock}

\end{proof}

\begin{remark}
\begin{enumerate}[label=(\roman*), ref=(\roman*), leftmargin=*, labelindent=0pt, itemsep=0pt, topsep=0.25em]
  \item Recall that $E(0)=\|\varphi-\bar{\phi}_L\|_{L^2(\R)}^2$. Since $\bar{\phi}_L-\bar{\phi}$ is supported in $[-L,L]$ for $L\geq 1$, 
  we have \begin{align*}
	E(0)&\leq 2\|\varphi-\bar{\phi}\|_{L^2(\R)}^2+2\|\bar{\phi}-\bar{\phi}_L\|_{L^2(\R)}^2\\
	&\leq 2\|\varphi-\bar{\phi}\|_{L^2(\R)}^2+16L\|\bar{\phi}\|_{L^\infty(\R)}^2\\
	&\leq CL,
  \end{align*}
where $C>0$ is independent of $L$. Thus, $E(0)=O(L)$. Consequently, for every fixed $T>0$, the bootstrap function $K_L(t)$ 
is also at most of order $L$ uniformly for $t\in [0,T]$.
  \item This observation explains the restriction $p=3,4$ in the global argument.
  For $p=3$, the coefficient of the superlinear term contains a negative power of $L$
  and can therefore be made sufficiently small by choosing $L$ larger.
  For $p=4$, the corresponding coefficient is of order one and remains uniformly bounded,
  which is just sufficient to close the bootstrap estimate.
  For $p>4$, however, the present $L^2$ energy argument fails to provide a uniform bound.

  \item Since the $H^1$ estimate in Lemma \ref{lemmaH1} is based essentially on the
  $L^2$ estimate established above, the continuation argument for global well-posedness
  is consequently restricted to $p=3,4$. For $p>4$, the local well-posedness result
  remains valid, but the present method does not yield the global a priori estimates
  required to extend the local solution globally.
\end{enumerate}
\end{remark}

\renewcommand{\theequation}{\thesection.\arabic{equation}}
\setcounter{equation}{0}
\section{Stability}
	\label{sec:stability}
	
We next examine the stability of traveling
waves governed by the gKdV-Burgers
equation \eqref{eq}. A fundamental feature of the traveling-wave profile is its translation
	invariance. More precisely, if $\widetilde{\phi}(\xi)$ is a traveling-wave
	profile, then, for any constant $a\in\mathbb{R}$,
	$
	\widetilde{\phi}(\xi-a)
	$
	represents the same shock profile translated in space. 
	Consequently, asymptotic stability with respect to a fixed, unshifted profile is impossible in general. 
	Indeed, for any $a\neq 0$, the translated traveling wave $\widetilde{\phi}(x-\sigma t -a)$ is itself an exact solution of \eqref{eq}, whereas
	$$\|\widetilde{\phi}(\cdot-\sigma t -a)-\widetilde{\phi}(\cdot-\sigma t)\|_{L^2(\R)}=\|\widetilde{\phi}(\cdot -a)-\widetilde{\phi}\|_{L^2(\R)}>0$$
	for all $t\geq 0$. Thus, even an exact shock differing only by a spatial translation does not converge to the fixed profile. 
	Hence stability must be formulated modulo translations, rather than with respect to a single representative of the traveling-wave family.

	To eliminate this translation-induced ambiguity, a time-dependent shift
	$Y(t)$ is introduced, and the solution is compared with the shifted profile
	\[
	\widetilde{\phi}(x-\sigma t-Y(t)).
	\]
	The function $Y(t)$ describes the displacement of the shock relative to its
	unperturbed position and will be specified in the corresponding contraction
	estimates.
	
To investigate the stability of the shock waves,
we derive an energy estimate measuring the
$L^2$-distance between the solution and a
time-dependently shifted shock profile.
The key evolution identity underlying the
subsequent analysis is stated in the proposition below.
	
	\begin{proposition}\label{pro:dt} Let \(\tilde{\phi}(\xi)\) be a monotone shock profile satisfying \eqref{boundcondition} and \eqref{eq:shockwave}. Let $\phi$ be a solution to \eqref{eq} such that
		$\phi-\bar{\phi}\in S_T$, and let
		$Y:[0,T]\to\mathbb{R}$ be Lipschitz continuous.
		Then, one obtains
		\begin{equation}\label{eq:dt}
			\begin{aligned}
				\frac{1}{2}\frac{d}{dt}  \int_{\mathbb{R}} (\phi^Y - \tilde{\phi})^2 \, d\xi
				&= - \epsilon \int_{\mathbb{R}} \left( (\phi^Y - \tilde{\phi})_\xi \right)^2 \, d\xi+\dot{Y}(t) \int_{\mathbb{R}} (\phi^Y - \tilde{\phi}) \tilde{\phi}' \, d\xi 
				 \\
				&- \int_{\mathbb{R}} \left( \frac{(\phi^Y)^p - \tilde{\phi}^p}{p} - \tilde{\phi}^{p-1} (\phi^Y - \tilde{\phi})\right) \tilde{\phi}'  d\xi
				,
			\end{aligned}
		\end{equation}
		where $\phi^Y(t, x) = \phi(t, x + Y(t)).$
	\end{proposition}
	
	\begin{proof}
		A direct computation of the time derivative of
		$\|\phi^Y-\tilde{\phi}\|_{L^2}^2$ gives
		$$\frac{d}{dt} \frac{1}{2} \int_{\mathbb{R}} (\phi^Y - \tilde{\phi})^2 \, d\xi
		=
		\frac{d}{dt} \frac{1}{2} \int_{\mathbb{R}} (\phi - \tilde{\phi}^{-Y})^2 \, d\xi. $$
		Set \( w = \phi - \tilde{\phi}^{-Y} \). 
		From  Eqs.~\eqref{eq} and \eqref{eq:shockwave}, one has
		\[\begin{aligned}
			w_t &= \phi_t - (\tilde{\phi}^{-Y})_t \\
			&= -\left(\frac{(w + \tilde{\phi}^{-Y})^p}{p}\right)_x + \epsilon (w + \tilde{\phi}^{-Y})_{xx} - \delta (w + \tilde{\phi}^{-Y})_{xxx}  + (\dot{Y}(t) + \sigma)(\tilde{\phi'})^{-Y} \\
			&=- \left( \frac{(w + \tilde{\phi}^{-Y})^p}{p} \right)_\xi 
			+ \epsilon w_{\xi\xi} - \delta w_{\xi\xi\xi} 
			+ \left( \frac{(\tilde{\phi}^{-Y})^p}{p} \right)_\xi 
			+ \dot{Y}(t)(\tilde{\phi'})^{-Y}. 
		\end{aligned}
		\]
		Then, it follows that 
		\[\begin{aligned}
			\frac{d}{dt} \frac{1}{2} \int_{\mathbb{R}} w^2 \, d\xi
			&= \int_{\mathbb{R}} w w_t \, d\xi \\
			&= \dot{Y}(t) \int_{\mathbb{R}} w (\tilde{\phi'})^{-Y} \, d\xi
			- \int_{\mathbb{R}} w\left( \frac{(w + \tilde{\phi}^{-Y})^p - (\tilde{\phi}^{-Y})^p}{p} \right)_\xi \, d\xi \\
			&+ \epsilon \int_{\mathbb{R}} w w_{\xi\xi} \, d\xi
			- \delta \int_{\mathbb{R}} w w_{\xi\xi\xi} \, d\xi \\
			&=\dot{Y}(t) \int_{\mathbb{R}} w (\tilde{\phi'})^{-Y} \, d\xi
			+ \int_{\mathbb{R}} w_\xi  \frac{(w + \tilde{\phi}^{-Y})^p - (\tilde{\phi}^{-Y})^p}{p}  \, d\xi - \epsilon \int_{\mathbb{R}}  w_{\xi} ^2 \, d\xi. 
		\end{aligned}	\]
		For convenience, set 
		\[
		I_1=\dot{Y}(t) \int_{\mathbb{R}} w (\tilde{\phi'})^{-Y} \, d\xi,\qquad
		I_2=\int_{\mathbb{R}} w_{\xi}\frac{\left(w+\tilde{\phi}^{-Y}\right)^{p}-\left(\tilde{\phi}^{-Y}\right)^{p}}{p} d\xi, \qquad
		I_3=\epsilon \int_{\mathbb{R}}  w_{\xi} ^2 \, d\xi, 
		\] and obtain
		\begin{align*}
			I_1&=\dot{Y}(t) \int_{\mathbb{R}} (\phi - \tilde{\phi}^{-Y}) (\tilde{\phi'})^{-Y} \, d\xi
			=\dot{Y}(t) \int_{\mathbb{R}} (\phi^{Y} - \tilde{\phi}) (\tilde{\phi'}) \, d(\xi+Y(t))
			=\dot{Y}(t) \int_{\mathbb{R}} (\phi^Y - \tilde{\phi}) \tilde{\phi}' \, d\xi , \\
			I_3&=\epsilon \int_{\mathbb{R}} ( (\phi - \tilde{\phi}^{-Y})_{\xi} )^2 \, d\xi
			=\epsilon \int_{\mathbb{R}} \left( (\phi^Y - \tilde{\phi})_\xi \right)^2 \, d\xi,  
		\end{align*}
			\begin{align*}
				I_2 &= \int_{\mathbb{R}} \frac{\phi^{p}-\left(\tilde{\phi}^{-Y}\right)^{p}}{p} \phi_{\xi} \, d\xi
				-\int_{\mathbb{R}} \frac{\phi^{p}-\left(\tilde{\phi}^{-Y}\right)^{p}}{p} (\tilde{\phi'})^{-Y} \, d\xi \\
				&= -\int_{\mathbb{R}} \phi\left(\phi^{p-1}\phi_{\xi}-\left(\tilde{\phi}^{-Y}\right)^{p-1}(\tilde{\phi'})^{-Y}\right) d\xi
				-\int_{\mathbb{R}} \frac{\phi^{p}-\left(\tilde{\phi}^{-Y}\right)^{p}}{p} (\tilde{\phi'})^{-Y} \, d\xi \\
				&= \frac{1}{p+1}\left[ \left(\tilde{\phi}^{-Y}\right)^{p+1}-\phi^{p+1} \right]\bigg|_{-\infty}^{+\infty}\\
				&\quad -\int_{\mathbb{R}} \left( \frac{\left(w+\tilde{\phi}^{-Y}\right)^{p}-\left(\tilde{\phi}^{-Y}\right)^{p}}{p} - \left(\tilde{\phi}^{-Y}\right)^{p-1} w \right) (\tilde{\phi'})^{-Y} \, d\xi \\
				&= -\int_{\mathbb{R}} \left( \frac{\left(w+\tilde{\phi}^{-Y}\right)^{p}-\left(\tilde{\phi}^{-Y}\right)^{p}}{p} - \left(\tilde{\phi}^{-Y}\right)^{p-1} w \right) (\tilde{\phi'})^{-Y}\, d\xi\\
				&=- \int_{\mathbb{R}} \left( \frac{(\phi^Y)^p - \tilde{\phi}^p}{p} - \tilde{\phi}^{p-1} (\phi^Y - \tilde{\phi})\right) \tilde{\phi}'  d\xi. 
			\end{align*}
		Thus, \eqref{eq:dt} follows.
		
	\end{proof}

	To apply the energy identity in Proposition~\ref{pro:dt}
	to the time-dependent shift used in the
	contraction estimates, we first establish
	the well-posedness of the shift equation.
In the following lemma, we verify the uniqueness
of the shift and its local Lipschitz continuity
with respect to time.
	\begin{lemma}[Well-posedness of the shift]
		\label{lem:shift-wellposed}
	Let $\widetilde{\phi}$ be a strictly decreasing shock profile connecting
	$\phi_{-\infty}$ and $\phi_{+\infty}$, where $\phi_{-\infty}>\phi_{+\infty}$,
	with $\widetilde{\phi}'\in L^1(\mathbb R)\cap L^2(\mathbb R)$.
		Suppose that $\phi$ is a solution of Eq.~\eqref{eq}
		on $[0,T_{\max})$ satisfying
		$\phi-\bar{\phi}\in S_T$ for every $T<T_{\max}$.
		Then the shift equation
		\begin{equation}
			\label{eq:shift-ode}
			\left\{
			\begin{aligned}
				\dot Y(t)
				={}&
				\frac{c}{\phi_{-\infty}-\phi_{+\infty}}
				\int_{\mathbb R}
				\bigl[\phi(t,x)
				-\widetilde{\phi}(x-\sigma t-Y(t))\bigr]
				\widetilde{\phi}'(x-\sigma t-Y(t))\,dx,\\
				Y(0)={}&0,
			\end{aligned}
			\right.
		\end{equation}
		where $c=c_1$ for $p=3$ and $c=c_2$ for $p=4$,
		admits a unique solution
		$Y\in C^1([0,T_{\max}))$.
		In particular, $Y$ is Lipschitz continuous
		on every compact subinterval of $[0,T_{\max})$.
	\end{lemma}
	
	\begin{proof}
	We apply the Cauchy-Lipschitz theorem.	Denote the right-hand side of
		\eqref{eq:shift-ode} by $F(t,Y)$.
		Setting $\xi=x-\sigma t$, we have
		\[
		\begin{aligned}
			F(t,Y)
			={}&
			\frac{c}{\phi_{-\infty}-\phi_{+\infty}}
			\int_{\mathbb R}
			\bigl[\phi(t,\xi+\sigma t)
			-\widetilde{\phi}(\xi-Y)\bigr]
			\widetilde{\phi}'(\xi-Y)\,d\xi.
		\end{aligned}
		\]
		By the translation invariance of
		\(
		\int_{\mathbb R}
		\widetilde{\phi}(\xi-Y)
		\widetilde{\phi}'(\xi-Y)\,d\xi,
		\)
		an application of integration by parts and
		the Cauchy-Schwarz inequality leads to
		\[
		\begin{aligned}
			|F(t,Y_1)-F(t,Y_2)|&=
			\frac{\left| c\right| }{\phi_{-\infty}-\phi_{+\infty}}
			\left| \int_{\mathbb{R}}
			\phi_x(t,\xi+\sigma t)
			\bigl[
			\widetilde{\phi}(\xi-Y_1)
			-\widetilde{\phi}(\xi-Y_2)
			\bigr]\,d\xi\right| \\
				&\leq
			\frac{|c|}
			{\phi_{-\infty}-\phi_{+\infty}}
			\|\phi_x(t)\|_{L^2(\R)}
			\|\widetilde{\phi}(\cdot-Y_1)
			-\widetilde{\phi}(\cdot-Y_2)\|_{L^2(\R)}\\
			&\quad\leq
			\frac{|c|}{\phi_{-\infty}-\phi_{+\infty}}
			\|\phi_x(t)\|_{L^2(\R)}
			\|\widetilde{\phi}'\|_{L^2(\R)}
			|Y_1-Y_2|.
		\end{aligned}
		\]
	Since $\phi-\bar{\phi}\in S_T$, the
	one-dimensional Sobolev embedding yields
	\[
	\sup_{0\leq t\leq T}
	\left(
	\|\phi(t)\|_{L^\infty(\R)}
	+\|\phi_x(t)\|_{L^2(\R)}
	\right)<\infty
	\]
	for every $T<T_{\max}$.
	Hence, $F$ is continuous in time and
	uniformly Lipschitz continuous with respect
	to $Y$ on $[0,T]$. By the Cauchy-Lipschitz
	theorem, the shift equation admits a
	unique local solution.
	Moreover, using
	\[
	\|\widetilde{\phi}'\|_{L^1(\R)}
	=\phi_{-\infty}-\phi_{+\infty},
	\]
	we obtain
	\[
	|F(t,Y)|
	\leq |c|
	\left(
	\|\phi(t)\|_{L^\infty(\R)}
	+\|\widetilde{\phi}\|_{L^\infty(\R)}
	\right)
	\leq M_T,
	\]
	where $M_T>0$ is independent of $Y$.
	Consequently,
	\[
	|Y(t)|\leq M_Tt,
	\]
	which implies that blow-up cannot occur in finite time.
	Consequently, the solution constructed locally admits
	a unique continuation to $[0,T_{\max})$ and remains
	Lipschitz continuous on each compact time interval.
	\end{proof}

	In view of the identity derived in
	Proposition~\ref{pro:dt}, the stability problem
	amounts to estimating the nonlinear flux
	contribution in terms of the viscous dissipation.
For this purpose, it is necessary to derive
quantitative pointwise bounds for the gradient
of the monotone traveling wave, which will be
obtained in the subsequent result.
	
	\begin{proposition}[Derivative estimates for the monotone shock profile]\label{pro:profile-derivative-general-p}
		Assume that the mapping \(a \mapsto a^{p-1}\) is strictly increasing on the interval \([\phi_{+\infty}, \phi_{-\infty}]\). 
		Let \(\tilde{\phi} \in C^3(\mathbb{R})\) be the strictly decreasing traveling-wave profile of  Eq.~\eqref{eq}. Assume further that \begin{equation*}
			\label{eq:discriminant-condition}
			0<4\delta\bigl(\phi_{-\infty}^{p-1}-\sigma\bigr)\leq\epsilon^2.
		\end{equation*} Then, for every \(\xi\in\mathbb R\), the following inequality holds:
		\begin{equation}\label{eq:profile-derivative-estimate}
			\begin{aligned}
				\kappa_1
				\left[
				\frac{\phi_{-\infty}^p-\tilde{\phi}^p(\xi)}{p}
				-\sigma\bigl(\phi_{-\infty}-\tilde{\phi}(\xi)\bigr)
				\right]
				\le -\epsilon\tilde{\phi}'(\xi)
				\le
				\kappa_2
				\left[
				\frac{\phi_{-\infty}^p-\tilde{\phi}^p(\xi)}{p}
				-\sigma\bigl(\phi_{-\infty}-\tilde{\phi}(\xi)\bigr)
				\right],
			\end{aligned}
		\end{equation}
		where
		\begin{equation*}\label{eq:kappa1-definition}
			\kappa_1=
			\frac{2}{
				1+\sqrt{1+
					\dfrac{4\delta\bigl(\sigma-\phi_{+\infty}^{p-1}\bigr)}{\epsilon^2}}},\qquad
			\kappa_2=
			\frac{2}{
				1+\sqrt{1-
					\dfrac{4\delta\bigl(\phi_{-\infty}^{p-1}-\sigma\bigr)}{\epsilon^2}}}.
		\end{equation*}
	\end{proposition}
	\begin{proof}
		Let
		\(Q(a):=\frac{\phi_{-\infty}^p-a^p}{p}-\sigma(\phi_{-\infty}-a)\), \(a\in[\phi_{+\infty},\phi_{-\infty}].\) 
		Therefore, it suffices to prove that
		\begin{equation*}
			\kappa_1Q(\tilde{\phi})\leq -\epsilon\tilde{\phi'} \leq 	\kappa_2Q(\tilde{\phi}). 
		\end{equation*}
		By the assumption, the nonlinear flux function \(f(\phi) = \frac{\phi^p}{p}\) is strictly convex on \([\phi_{+\infty}, \phi_{-\infty}]\), and it follows that
	\begin{equation*}\label{eq:characteristic-speed-order}
		\phi_{+\infty}^{p-1}<\sigma=\dfrac{\phi_{-\infty}^p-\phi_{+\infty}^p}{p(\phi_{-\infty}-\phi_{+\infty})}<\phi_{-\infty}^{p-1}, \qquad
		Q(\widetilde{\phi})>0.
	\end{equation*}
		The integrated traveling-wave equation is
		\[
		\epsilon\widetilde{\phi}'
		-\delta\widetilde{\phi}''
		=-Q(\widetilde{\phi}).
		\]
		For $\kappa>0$, define
		$
		Z_\kappa(\xi):=
		\epsilon\widetilde{\phi}'(\xi)
		+\kappa Q(\widetilde{\phi}(\xi)).
		$
	Differentiating the traveling wave equation above directly yields
		\begin{equation}\label{Z-kappa}
		Z_\kappa'-A_\kappa Z_\kappa
		=B_\kappa Q(\widetilde{\phi}),
	\end{equation}
		where
		\[
		A_\kappa:=
		\frac{\epsilon}{\delta}
		+\frac{\kappa}{\epsilon}
		Q'(\widetilde{\phi}),
		\qquad
		B_\kappa:=
		\frac{\epsilon}{\delta}(1-\kappa)
		-\frac{\kappa^2}{\epsilon}
		Q'(\widetilde{\phi}).
		\]
		Since
		\[
		Q'(a)=\sigma-a^{p-1},
		\qquad
		0<4\delta(\phi_{-\infty}^{p-1}-\sigma)
		\leq\epsilon^2,
		\]
		and $0<\kappa_1<1<\kappa_2\leq2$, it follows that
		\begin{equation}\label{A-kappa}
		\begin{aligned}
			A_\kappa
			\geq
			\frac{\epsilon}{\delta}
			-\frac{\kappa}{\epsilon}
			(\phi_{-\infty}^{p-1}-\sigma)
			\geq
			\frac{\epsilon}{\delta}
			\left(1-\frac{\kappa}{4}\right)
			\geq\frac{\epsilon}{2\delta}>0,
		\end{aligned}
	\end{equation}
		for $\kappa=\kappa_1,\kappa_2$.
		Multiplying \eqref{Z-kappa} by the integrating factor
		$
		\exp\left(-\int_\xi^s A_\kappa(\tau)\,d\tau\right)
		$
		and integrating from $\xi$ to $R$, we obtain
		\[
		\begin{aligned}
			Z_\kappa(\xi)
			={}
			e^{-\int_\xi^R A_\kappa(\tau)\,d\tau}
			Z_\kappa(R)
			-
			\int_\xi^R
			e^{-\int_\xi^s A_\kappa(\tau)\,d\tau}
			B_\kappa(s)Q(\widetilde{\phi}(s))\,ds.
		\end{aligned}
		\]
		Since $Z_\kappa(R)\to0$ as $R\to+\infty$,
		\eqref{A-kappa} allows us to pass to the limit 
		\begin{equation}\label{Z-kappa1}
			Z_\kappa(\xi)
			=
			-\int_\xi^{+\infty}
			e^{-\int_\xi^s A_\kappa(\tau)\,d\tau}
			B_\kappa(s)Q(\widetilde{\phi}(s))\,ds.
		\end{equation}
		By the definitions of $\kappa_1$ and $\kappa_2$,
		\[
		\frac{\epsilon}{\delta}(1-\kappa_1)
		=
		\frac{\kappa_1^2}{\epsilon}
		(\sigma-\phi_{+\infty}^{p-1}), \qquad
		\frac{\epsilon}{\delta}(1-\kappa_2)
		=
		-\frac{\kappa_2^2}{\epsilon}
		(\phi_{-\infty}^{p-1}-\sigma).
		\]
		Consequently,
		\[
			B_{\kappa_1}
			=
			\frac{\kappa_1^2}{\epsilon}
			\bigl(
			\widetilde{\phi}^{p-1}
			-\phi_{+\infty}^{p-1}
			\bigr)>0, \qquad
				B_{\kappa_2}
			=
			-\frac{\kappa_2^2}{\epsilon}
			\bigl(
			\phi_{-\infty}^{p-1}
			-\widetilde{\phi}^{p-1}
			\bigr)<0.
		\]
		Since the exponential factor and
		$Q(\widetilde{\phi})$ in \eqref{Z-kappa1} are positive,
		it follows that
		\[
		Z_{\kappa_1}(\xi)\leq0,
		\qquad
		Z_{\kappa_2}(\xi)\geq0.
		\]
		Thus, by the definition of $Z_\kappa$,
		\[
			\kappa_1Q(\widetilde{\phi}(\xi))
			\leq
			-\varepsilon\widetilde{\phi}'(\xi)
			\leq
			\kappa_2Q(\widetilde{\phi}(\xi)),
			\qquad \xi\in\mathbb R.
		\]
	
		The argument also applies when
		$\epsilon^2=
		4\delta(\phi_{-\infty}^{p-1}-\sigma)$,
		since then $\kappa_2=2$ and \eqref{Z-kappa}-\eqref{Z-kappa1}
		remain valid. 
	\end{proof}

	\begin{remark}
		When \(p\) is even, \(p-1\) is odd, so the mapping \(a \mapsto a^{p-1}\) is strictly increasing on \(\mathbb{R}\). If \(p\) is odd, one usually needs to restrict the end states to a convex interval of the flux; for instance, taking \(0 \le \phi_{+\infty} < \phi_{-\infty}\) suffices to guarantee the above monotonicity condition.
	\end{remark}
	
	\begin{corollary}\label{p=3-4,corollary}
		The assumptions of
		Proposition~\ref{pro:profile-derivative-general-p}
		lead to the explicit bounds stated below.
		When \(p = 3\), 
		\begin{equation}\label{eq:p=3-inequality}
			\lambda_1(\phi_{-\infty} -\tilde{\phi})(\tilde{\phi}-\phi_{+\infty})<-\epsilon \tilde{\phi}'<\lambda_2(\phi_{-\infty} -\tilde{\phi})(\tilde{\phi}-\phi_{+\infty}),  
		\end{equation}
		where $\lambda_1=\frac{2(\sqrt{2}-1)}{3}(\phi_{-\infty}+2\phi_{+\infty})$, $\lambda_2=\frac{2}{3}(2\phi_{-\infty}+\phi_{+\infty})$. 
		
		When \(p = 4\), factoring $Q(\tilde{\phi})$ gives
		\[
		\frac{\phi_{-\infty}^4-\widetilde\phi^4}{4}
		-\sigma(\phi_{-\infty}-\widetilde\phi)
		=
		\frac14
		(\phi_{-\infty}-\widetilde\phi)
		(\widetilde\phi-\phi_{+\infty})
		H(\widetilde\phi),
		\]
		where
		$
		H(z)
		=
		z^2+(\phi_{-\infty}+\phi_{+\infty})z
		+\phi_{-\infty}^2
		+\phi_{-\infty}\phi_{+\infty}
		+\phi_{+\infty}^2.
		$
	Define
	\[
	m_4:=\min_{z\in[\phi_{+\infty},\phi_{-\infty}]}H(z),
	\qquad
	M_4:=\max_{z\in[\phi_{+\infty},\phi_{-\infty}]}H(z),
	\]
	so that
	$0<m_4\le M_4$
	(since $\Delta<0$, we have $H(z)>0$ for all $z\in\mathbb R$).
	Hence,
		\begin{equation}\label{eq:p'=4-inequality}
		\lambda_3
		(\phi_{-\infty}-\widetilde\phi)
		(\widetilde\phi-\phi_{+\infty})
		\le
		-\epsilon\widetilde\phi'
		\le
		\lambda_4
		(\phi_{-\infty}-\widetilde\phi)
		(\widetilde\phi-\phi_{+\infty}),
	\end{equation}
		where
		$
		\lambda_3=\frac{\kappa_1m_4}{4}$ and
		$
		\lambda_4=\frac{\kappa_2M_4}{4}.
		$
		
		In particular, if $\phi_{+\infty}>0$, a more precise estimate can be obtained,
		\begin{equation}\label{eq:p=4-inequality}
			\lambda_3'(\phi_{-\infty} -\tilde{\phi})(\tilde{\phi}-\phi_{+\infty})<-\epsilon \tilde{\phi}'<\lambda_4'(\phi_{-\infty} -\tilde{\phi})(\tilde{\phi}-\phi_{+\infty}),  
		\end{equation}
		where $\lambda_3'=\frac{\sqrt{2}-1}{2}(\phi_{-\infty}^2+2\phi_{-\infty}\phi_{+\infty} +3\phi_{+\infty}^2)$, $\lambda_4'=\frac{1}{2}(3\phi_{-\infty}^2+2\phi_{-\infty}\phi_{+\infty} +\phi_{+\infty}^2)$. 
	\end{corollary}
	
	Proposition~\ref{pro:profile-derivative-general-p} and Corollary~\ref{p=3-4,corollary} provide the
	pointwise estimates for the shock profile required in the subsequent
	stability analysis. Since $\widetilde{\phi}$ is strictly decreasing, it holds that
	\[
	\|\widetilde{\phi}'\|_{L^1(\R)}
	=
	-\int_{\mathbb R}\widetilde{\phi}'(\xi)\,d\xi
	=
	\phi_{-\infty}-\phi_{+\infty}
	<\infty.
	\]
	Moreover, the upper profile estimates in Corollary~\ref{p=3-4,corollary} imply
	$
	\|\widetilde{\phi}'\|_{L^\infty(\R)}<\infty.
	$
	Hence, by H\"{o}lder's inequality, it follows that 
	\[
	\|\widetilde{\phi}'\|_{L^2(\R)}^2
	\le
	\|\widetilde{\phi}'\|_{L^\infty(\R)}
	\|\widetilde{\phi}'\|_{L^1(\R)}
	<\infty,
	\]
	and therefore
	\[
	\widetilde{\phi}'
	\in L^1(\mathbb R)\cap L^2(\mathbb R).
	\]
	
	These profile estimates, together with the energy identity in
	Proposition~\ref{pro:dt} and the Poincar\'e-type inequality after the change
	of variables $z=\widetilde{\phi}(\xi)$, are used below to establish
	the $L^2$-contraction estimates for $p=3$ and $p=4$. 
	We first state the following lemmas.

	\begin{lemma}
		\label{lem:bootstrap-closure}
	Let $\tilde{\phi}$ be a strictly decreasing shock profile and let $\phi$
	be the corresponding global solution. Let $Y:[0,\infty)\to\mathbb{R}$ be the corresponding Lipschitz shift satisfying $Y(0)=0$. Assume that, for every finite
	time interval $[0,T]$ on which the solution is defined,
		\[
		\sup_{0\le t\le T}
		\|\phi(t,\cdot)-\tilde{\phi}(\cdot-\sigma t-Y(t))\|_{L^\infty(\mathbb R)}
		\le \Gamma,
		\]
		the estimates
		\[
		\|\phi(t,\cdot)-\tilde{\phi}(\cdot-\sigma t-Y(t))\|_{L^2(\R)}^2
		+
		c_0\epsilon
		\int_0^t
		\|\partial_x[
		\phi(s,\cdot)-\tilde{\phi}(\cdot-\sigma s-Y(s))
		]\|_{L^2(\R)}^2\,ds
		\le
		\|\varphi-\tilde{\phi}\|_{L^2(\R)}^2
		\]
		and
		\begin{equation*}
			\begin{aligned}
				\frac{d}{dt}
				\|\partial_x[
				\phi(t,\cdot)-\tilde{\phi}(\cdot-\sigma t-Y(t))
				]\|_{L^2(\R)}^2
				&\le
				C_\Gamma
				(
				\|\partial_x[
				\phi(t,\cdot)-\tilde{\phi}(\cdot-\sigma t-Y(t))
				]\|_{L^2(\R)} \\
				&+
				\|\partial_x[
				\phi(t,\cdot)-\tilde{\phi}(\cdot-\sigma t-Y(t))
				]\|_{L^2(\R)}^2
				)
			\end{aligned}
		\end{equation*}
		hold, where $c_0,l>0$ and $C_\Gamma>0$ are independent of $T$. Then there exists $\varepsilon_0>0$ such that
		\[
		\|\varphi-\tilde{\phi}\|_{H^1(\mathbb R)}
		\le \varepsilon_0
		\]
		implies
		\[
		\sup_{t\ge0}
		\|\phi(t,\cdot)-\tilde{\phi}(\cdot-\sigma t-Y(t))\|_{L^\infty(\mathbb R)}
		<
		\Gamma.
		\]
	\end{lemma}
	
	\begin{proof}
		Set
		$
		w(t,x)
		:=
		\phi(t,x)-\tilde{\phi}(x-\sigma t-Y(t))$ and $
			D(t):=\|w_x(t)\|_{L^2(\mathbb R)}^2
		$. 
		Since $Y(0)=0$, 
		$
		w(0,x)=\varphi(x)-\tilde{\phi}(x).
		$
		By the Sobolev embedding $H^1(\mathbb{R}) \hookrightarrow L^\infty(\mathbb{R})$,
		\[
		\|w(0)\|_{L^\infty(\R)}^2
		\le
		\|w(0)\|_{H^1}^2.
		\]
		Hence, if $\varepsilon_0>0$ is sufficiently small,
		it follows that $
		\|w(0)\|_{L^\infty(\R)}
		<
		\frac{\Gamma}{2}.
		$
		Define
		\[
		T^*
		:=
		\sup
		\left\{
		T>0:
		\sup_{0\le t\le T}
		\|w(t)\|_{L^\infty(\mathbb R)}
		<
		\Gamma
		\right\}.
		\]
		Since
		$
		w\in C([0,\infty);H^1(\mathbb R))
		\hookrightarrow
		C([0,\infty);L^\infty(\mathbb R)),
		$
		we have $T^*>0$. Fix any $T<T^*$, and by the definition of $T^*$,
		\[
		\sup_{0\le t\le T}\|w(t)\|_{L^\infty(\R)}<\Gamma.
		\]
		Therefore, the assumed $L^2$ estimate gives
		\[
		\sup_{0\le t\le T}
		\|w(t)\|_{L^2(\R)}
		\le
		\|w(0)\|_{L^2(\R)}
		\le
		\varepsilon_0\quad \text{and} \quad
		\int_0^T D(t)\,dt
		\le
		\frac{\varepsilon_0^2}{c_0\epsilon}.
		\]
		
		We next prove that
		\[
		D(t)<1,
		\qquad
		0\le t\le T.
		\]
		Choose $\varepsilon_0$ so small that
	\[
    \varepsilon_0^2<\min\left\{\frac{1}{2},\frac{c_0\epsilon}{8C_\Gamma},\frac{\Gamma^2}{4l}\right\}.
	\]
		Then
		\[
		D(0)
		=
		\|w_x(0)\|_{L^2(\R)}^2
		\le
		\varepsilon_0^2
		<
		\frac12.
		\]
		Suppose that $D$ reaches the value $1$ on $[0,T]$.
		Then, by continuity, there exist $0\le s<t\le T$ such that
		$
		D(s)=\frac12$, 
		$
		D(t)=1,
		$
		and
		$
		\frac12\le D(\tau)\le1$,
		$
		\tau\in[s,t].
		$
		Using the assumed differential inequality,
		\[
		D'(\tau)
		\le
		C_\Gamma
		\bigl(D(\tau)^{1/2}+D(\tau)\bigr)
		\le
		2C_\Gamma.
		\]
		Hence, integrating the above inequality yields
		\[
		\frac12
		=
		D(t)-D(s)
		\le
		2C_\Gamma(t-s),
		\]
		so that
		$
		t-s
		\ge
		\frac{1}{4C_\Gamma}.
		$
		Therefore,
		\[
		\int_s^t D(\tau)\,d\tau
		\ge
		\frac12(t-s)
		\ge
		\frac{1}{8C_\Gamma}.
		\]
		On the other hand,
		\[
		\int_s^t D(\tau)\,d\tau
		\le
		\int_0^T D(\tau)\,d\tau
		\le
		\frac{\varepsilon_0^2}{c_0\epsilon}.
		\]
		This is impossible if 
		$
		\frac{\varepsilon_0^2}{c_0\epsilon}
		<
		\frac{1}{8C_\Gamma}.
		$
		Thus
		$
		D(t)<1$,
		$
		0\le t\le T.
		$
		
		Then, using the one-dimensional Gagliardo-Nirenberg inequality,
		it holds
		\[
		\|w(t)\|_{L^\infty(\R)}^2
		\le
		l\|w(t)\|_{L^2(\R)}\|w_x(t)\|_{L^2(\R)}
		\le
		l\varepsilon_0,
		\qquad
		0\le t\le T.
		\]
		Choose $\varepsilon_0$ further so that
		$
		l\varepsilon_0
		<
		\frac{\Gamma^2}{4}.
		$
		Then
		\[
		\sup_{0\le t\le T}
		\|w(t)\|_{L^\infty(\R)}
		<
		\frac{\Gamma}{2}.
		\]
		Since the constants above are independent of $T<T^*$, letting
		$T\rightarrow T^*$ yields
		\[
		\sup_{0\le t<T^*}
		\|w(t)\|_{L^\infty(\R)}
		\le
		\frac{\Gamma}{2}.
		\]
		If $T^*<\infty$, continuity in $L^\infty(\mathbb R)$ implies
		$
		\|w(T^*)\|_{L^\infty(\R)}
		\le
		\frac{\Gamma}{2}.
		$
		There is therefore a positive number $\tau$
		for which
		\[
		\|w(t)\|_{L^\infty(\R)}<\Gamma
		\quad \text{for all }
		t\in[T^*,T^*+\tau].
		\]
		Such a continuation is incompatible with
		the definition of $T^*$. Hence $T^*$ cannot
		be finite, which implies
		\[
		T^*=\infty.
		\]
		As a result, the perturbation satisfies
		\[
		\sup_{t\ge0}
		\|w(t)\|_{L^\infty(\R)}<\Gamma.
		\]
		This concludes the argument.
		
	\end{proof}
	
	\begin{lemma}[Derivative estimate]
		\label{lem:derivative-estimate}
		Let $p\in\{3,4\}$, and let $\widetilde{\phi}$ be a strictly decreasing
		shock profile satisfying
		$
		\widetilde{\phi}'\in L^1(\mathbb{R})\cap L^2(\mathbb{R}).
		$
		Assume that, for some $M,\Gamma>0$,
		\[
		\sup_{0\le t\le T}
		\|\phi(t,\cdot)-\widetilde{\phi}(\cdot-\sigma t-Y(t))\|_{L^2(\R)}
		\le M,
		\]
		\[
		\sup_{0\le t\le T}
		\|\phi(t,\cdot)-\widetilde{\phi}(\cdot-\sigma t-Y(t))\|_{L^\infty(\R)}
		\le \Gamma, 
		\]
		and
	\[
	|\dot Y(t)|\leq
	c\|\phi(t,\cdot)-\widetilde{\phi}(\cdot-\sigma t-Y(t))\|_{L^\infty(\R)}.  
	\]
		Then
		\[
		\begin{aligned}
			&\frac{d}{dt}
			\left\|
			\partial_x\bigl[
			\phi(t,\cdot)-\widetilde{\phi}(\cdot-\sigma t-Y(t))
			\bigr]
			\right\|_{L^2(\R)}^2
			+\epsilon
			\left\|
			\partial_x^2\bigl[
			\phi(t,\cdot)-\widetilde{\phi}(\cdot-\sigma t-Y(t))
			\bigr]
			\right\|_{L^2(\R)}^2
			\\
			&\le
			C\left(
			\left\|
			\partial_x\bigl[
			\phi(t,\cdot)-\widetilde{\phi}(\cdot-\sigma t-Y(t))
			\bigr]
			\right\|_{L^2(\R)}
			+
			\left\|
			\partial_x\bigl[
			\phi(t,\cdot)-\widetilde{\phi}(\cdot-\sigma t-Y(t))
			\bigr]
			\right\|_{L^2(\R)}^2
			\right),
		\end{aligned}
		\]
		where $C>0$ is independent of $T$.
	\end{lemma}
	
	\begin{proof}
		For convenience, write
		\[
		\widetilde{\phi}^{-Y}
		=\widetilde{\phi}(x-\sigma t-Y(t)),
		\qquad
		w=\phi-\widetilde{\phi}^{-Y},
		\]
		and let
		\[
		U_*=\max\left\{
		|\phi_{-\infty}|,|\phi_{+\infty}|
		\right\}.
		\]
		The evolution identity established in
		Proposition~\ref{pro:dt} then yields
		\[
		w_t
		=
		-\partial_x\left[
		\frac{(w+\widetilde{\phi}^{-Y})^p
			-(\widetilde{\phi}^{-Y})^p}{p}
		\right]
		+\epsilon w_{xx}
		-\delta w_{xxx}
		+\dot Y(t)(\widetilde{\phi}^{-Y})_x.
		\]
		Multiplying the equation by $-w_{xx}$ and integrating over $\mathbb{R}$ gives
		\[
		\begin{aligned}
			\frac12\frac{d}{dt}\|w_x\|_{L^2(\R)}^2
			+\epsilon\|w_{xx}\|_{L^2(\R)}^2
			=
			\int_{\mathbb R}w_{xx}
			\partial_x\left[
			\frac{(w+\widetilde{\phi}^{-Y})^p
				-(\widetilde{\phi}^{-Y})^p}{p}
			\right]dx
			-\dot Y(t)
			\int_{\mathbb R}
			w_{xx}(\widetilde{\phi}^{-Y})_x\,dx .
		\end{aligned}
		\]
	By the assumption on \([0,T]\), one obtains
		\[
		\|w(t)\|_{L^2(\R)}\leq M, \qquad \|w(t)\|_{L^\infty(\R)}\leq\Gamma.
		\]
			Thus, $|\phi(t,x)|\leq U_*+\Gamma$ and
		$|\widetilde{\phi}^{-Y}|\leq U_*$.
		Since $p\in\{3,4\}$, the mean value theorem yields
		\[
		\left|
		\phi^{p-1}-(\widetilde{\phi}^{-Y})^{p-1}
		\right|
		\leq (p-1)(U_*+\Gamma)^{p-2}|w|.
		\]
		Consequently,
		\[
		\begin{aligned}
			\left\|
			\partial_x\left[
			\frac{(w+\widetilde{\phi}^{-Y})^p
				-(\widetilde{\phi}^{-Y})^p}{p}
			\right]
			\right\|_{L^2(\mathbb{R})}^2
			&=
			\left\|
			\phi^{p-1}w_x+
			\left(\phi^{p-1}
			-(\widetilde{\phi}^{-Y})^{p-1}\right)
			(\widetilde{\phi}^{-Y})_x
			\right\|_{L^2(\mathbb{R})}^2\\
			&\leq
			2(U_*+\Gamma)^{2(p-1)}
			\|w_x\|_{L^2(\mathbb{R})}^2\\
			&\quad+
			2(p-1)^2(U_*+\Gamma)^{2(p-2)}
			\|\widetilde{\phi}'\|_{L^2(\mathbb{R})}^2
			\|w\|_{L^\infty(\mathbb{R})}^2\\
			&\leq
			C_{U_*,\Gamma,p,\|\widetilde{\phi}'\|_{L^2(\mathbb{R})}}
			\left(
			\|w_x\|_{L^2(\mathbb{R})}^2+
			\|w\|_{L^\infty(\mathbb{R})}^2
			\right).
		\end{aligned}
		\]
		The one-dimensional Gagliardo-Nirenberg
		interpolation inequality also ensures the
		existence of a constant $l>0$ satisfying
		\[
		\|w\|_{L^\infty(\R)}^2
		\le
		lM\|w_x\|_{L^2(\R)}.
		\]
		Hence, by applying the Cauchy–Schwarz and Young inequalities, it follows that
		\[
		\left|
		\int_{\mathbb R}w_{xx}
		\partial_x\left[
		\frac{(w+\widetilde{\phi}^{-Y})^p
			-(\widetilde{\phi}^{-Y})^p}{p}
		\right]dx
		\right|
		\le
		\frac{\epsilon}{4}\|w_{xx}\|_{L^2(\R)}^2
		+
		C_{U_*,\Gamma,p,l,M,\|\tilde{\phi'}\|_{L^2(\R)}}\left(
		\|w_x\|_{L^2(\R)}
		+\|w_x\|_{L^2(\R)}^2
		\right).
		\]
		Similarly,
		\[
		|\dot Y(t)|
		\le
			\dfrac{|c|}{\phi_{-\infty}-\phi_{+\infty}}\|w\|_{L^\infty(\R)}\|\tilde{\phi'}\|_{L^1(\R)}
		\le
		|c|\|w_x\|_{L^2(\R)}^{1/2},
		\]
		so that
		\[
		\left|
		\dot Y(t)
		\int_{\mathbb R}
		w_{xx}(\widetilde{\phi}^{-Y})_x\,dx
		\right|
		\le
		\frac{\epsilon}{4}\|w_{xx}\|_{L^2(\R)}^2
		+
		|c|\|w_x\|_{L^2(\R)}.
		\]
		Combining the above estimates yields
		\[
		\frac{d}{dt}\|w_x\|_{L^2(\R)}^2
		+
		\epsilon\|w_{xx}\|_{L^2(\R)}^2
		\le
			C_{U_*,\Gamma,p,l,M,c,\|\tilde{\phi'}\|_{L^2(\R)}}\left(
		\|w_x\|_{L^2(\R)}
		+
		\|w_x\|_{L^2(\R)}^2
		\right),
		\]
		which proves the result.
	\end{proof}
	
\begin{lemma}[Poincar\'e-type inequality]\cite{chen20262,MR4195742}\label{Weighted Poincaré inequality}
	For any function \(f:[a,b]\to\mathbb{R}\) satisfying
	$
	\int_a^b (y-a)(b-y)|f'|^2\,dy<\infty,
	$
	one has
	\[
	\int_a^b f^2\,dy
	\le
	\frac12\int_a^b (y-a)(b-y)|f'|^2\,dy
	+
	\frac{1}{b-a}
	\left(\int_a^b f\,dy\right)^2.
	\]
\end{lemma}
	
	Next, the proofs of Theorems \ref{thm:main} and \ref{the:p'=4} are presented below.

\begin{proof} ({\bf Proof of Theorem \ref{thm:main}})
		Set
		$
		w(t,x)
		:=
		\phi(t,x)
		-\widetilde{\phi}(x-\sigma t-Y(t)).
		$
		Let $T>0$ and assume temporarily that
		\[
		\sup_{0\le t\le T}
		\|w(t)\|_{L^{\infty}(\mathbb{R})}
		\le\gamma.
		\]
		By Proposition~\ref{pro:dt}, for $p=3$,
		\[
		\frac12\frac{d}{dt}
		\int_{\mathbb{R}}w^{2}\,d\xi
		=
		\dot{Y}(t)
		\int_{\mathbb{R}}w\widetilde{\phi}'\,d\xi
		-
		\int_{\mathbb{R}}
		\left(
		\frac{w}{3}+\widetilde{\phi}
		\right)
		w^{2}\widetilde{\phi}'\,d\xi
		-
		\epsilon\int_{\mathbb{R}}w_{\xi}^{2}\,d\xi.
		\]
		Introduce the following variables
		$
		z=\widetilde{\phi}(\xi)$, 
		$
		dz=\widetilde{\phi}'(\xi)\,d\xi$, and 
		$
			v(z)=w\circ\widetilde{\phi}^{-1}(z).
		$
		Then
			\begin{align*}
			\dot{Y}(t)\int_{\mathbb{R}} w\,\tilde{\phi}'\,d\xi
			&= \frac{c_1}{\phi_{-\infty} - \phi_{+\infty}}\left( \int_{\phi_{+\infty}}^{\phi_{-\infty}} v(z)\,dz\right) ^2, \\
			- \int_{\mathbb{R}} \left( \frac{w}{3}+\tilde{\phi}\right) w^2 \tilde{\phi}'\,d\xi
			&=\int_{\phi_{+\infty}}^{\phi_{-\infty}} \left( \frac{v}{3}+z\right)v(z)^2\,dz
			\leq \frac{
				\gamma+3\phi_{-\infty}}{3}\int_{\phi_{+\infty}}^{\phi_{-\infty}} v(z)^2\,dz, \\
			-\epsilon \int_{\mathbb{R}} w_\xi^2 \, d\xi
			&= -\int_{\phi_{+\infty}}^{\phi_{-\infty}} (-\epsilon \tilde{\phi}') v_z^2 \, dz.
		\end{align*}
		By inequality \eqref{eq:p=3-inequality} in Corollary~\ref{p=3-4,corollary}, it follows that
		\[
		-\epsilon\int_{\mathbb{R}}w_{\xi}^{2}\,d\xi
		\le
		-\lambda_{1}
		\int_{\phi_{+\infty}}^{\phi_{-\infty}}
		(\phi_{-\infty}-z)(z-\phi_{+\infty})
		v_{z}^{2}\,dz.
		\]
		Using the Poincar\'e-type inequality in Lemma \ref{Weighted Poincaré inequality},  
		it holds that
		\begin{align*}
			\frac{d}{dt}\frac{1}{2}\int_{\mathbb{R}} w^2 \, d\xi
			&\le -\frac{\gamma+3\phi_{-\infty}}{3(\phi_{-\infty} - \phi_{+\infty})}\left(\int_{\phi_{+\infty}}^{\phi_{-\infty}} v\,dz\right)^2
			+ \frac{\gamma+3\phi_{-\infty}}{6}\int_{\phi_{+\infty}}^{\phi_{-\infty}} (\phi_{-\infty} - z)(z - \phi_{+\infty}) v_z^2 \, dz \\
			&
			\quad+ \frac{\gamma+3\phi_{-\infty}}{3(\phi_{-\infty} - \phi_{+\infty})}\left(\int_{\phi_{+\infty}}^{\phi_{-\infty}} v\,dz\right)^2 
			- \lambda_1 \int_{\phi_{+\infty}}^{\phi_{-\infty}}  (\phi_{-\infty} - z)(z - \phi_{+\infty}) v_z^2 \, dz \\
			&=  -(\lambda_1-\dfrac{\gamma+3\phi_{-\infty}}{6})
			\int_{\phi_{+\infty}}^{\phi_{-\infty}}  (\phi_{-\infty} - z)(z - \phi_{+\infty}) v_z^2 \, dz \le 0.
		\end{align*}
		Using the upper profile estimate \eqref{eq:p=3-inequality},
		\[
		\epsilon\|w_x(t)\|_{L^{2}}^{2}
		\le
		\lambda_{2}
		\int_{\phi_{+\infty}}^{\phi_{-\infty}}
		(\phi_{-\infty}-z)(z-\phi_{+\infty})
		v_{z}^{2}\,dz,
		\]
		we conclude that
		\[
		\frac12\frac{d}{dt}\|w(t)\|_{L^{2}}^{2}
		+
		A_{1}^{*}\epsilon\|w_x(t)\|_{L^{2}}^{2}
		\le0,
		\]
		where
		$
		A_{1}^{*}
		=
		\frac{8(\sqrt{2}-1)\phi_{+\infty}+(4\sqrt{2}-7)\phi_{-\infty}-\gamma}{6\lambda_{2}}>0.
		$
		Integrating in time gives
		\[
		\|w(t)\|_{L^{2}}^{2}
		+
		2A_{1}^{*}\epsilon
		\int_{0}^{t}\|w_x(s)\|_{L^{2}}^{2}\,ds
		\le
		\|w(0)\|_{L^{2}}^{2},
		\qquad
		0\le t\le T.
		\]
		In particular,
		\[
		\sup_{0\le t\le T}\|w(t)\|_{L^{2}}
		\le
		\|w(0)\|_{L^{2}}.
		\]
		
		Next, since 
		$
		\|\widetilde{\phi}'\|_{L^{1}}
		=
		\phi_{-\infty}-\phi_{+\infty}
		$, the definition of $Y(t)$ gives
		\[
		\begin{aligned}
			|\dot{Y}(t)|
			&\le
			\frac{|c_{1}|}{\phi_{-\infty}-\phi_{+\infty}}
			\|w(t)\|_{L^{\infty}}
			\|\widetilde{\phi}'\|_{L^{1}}
			\le
			|c_{1}|\|w(t)\|_{L^{\infty}}.
		\end{aligned}
		\]
		Therefore, under the bootstrap assumption,
		Lemma~\ref{lem:derivative-estimate} yields
		\[
		\frac{d}{dt}\|w_x(t)\|_{L^{2}}^{2}
		\le
		C_{\gamma}
		\left(
		\|w_x(t)\|_{L^{2}}
		+
		\|w_x(t)\|_{L^{2}}^{2}
		\right),
		\]
		where $C_{\gamma}>0$ is independent of $T$.
		Lemma~\ref{lem:bootstrap-closure} is therefore
		applicable with $\Gamma=\gamma$ and $c_0=2A_1^*$.
		In particular, there is an $\varepsilon_0>0$
		such that
		\(
		\|\varphi-\widetilde{\phi}\|_{H^1(\mathbb{R})}
		<\varepsilon_0
		\)
		yields
		\[
		\sup_{t\ge0}\|w(t)\|_{L^{\infty}(\mathbb{R})}
		<\gamma.
		\]
		
		Therefore the conditional $L^{2}$-contraction estimate obtained above
		holds on every finite time interval and hence, for every $t\ge0$,
		\[
		\|w(t)\|_{L^{2}}^{2}
		+
		2A_{1}^{*}\epsilon
		\int_{0}^{t}\|w_x(s)\|_{L^{2}}^{2}\,ds
		\le
		\|w(0)\|_{L^{2}}^{2}.
		\]

		\end{proof}
\begin{proof} ({\bf 	Proof of Theorem \ref{the:p'=4}})
			Set
		$
		w(t,x)
		=
		\phi(t,x)
		-
		\widetilde{\phi}(x-\sigma t-Y(t)).
		$
		On a finite interval $[0,T]$, assume temporarily that
		\[
		\sup_{0\le t\le T}\|w(t)\|_{L^\infty(\R)}\le\gamma.
		\]
		Using the same change of variables
		$
		z=\widetilde{\phi}(\xi)
		$
		as in Theorem~\ref{thm:main}, we have $|v|\le\gamma$. Since
		$|z|\le U_*$,
		\[
		\frac14v^4+zv^3+\frac32z^2v^2
		\le
		B_\gamma v^2,
		\qquad
		B_\gamma
		:=
		\frac{\gamma^2}{4}
		+\gamma U_*
		+\frac32U_*^2.
		\]
	With the choice $c_2=-B_\gamma$, using Proposition~\ref{pro:dt}, Lemma \ref{Weighted Poincaré inequality}, and the two-sided profile estimate \eqref{eq:p'=4-inequality} in Corollary \ref{p=3-4,corollary} gives
		\[
		\frac12\frac{d}{dt}\|w(t)\|_{L^2(\R)}^2
		\le
		-\left(
		\lambda_3-\frac{B_\gamma}{2}
		\right)
		\int_{\phi_{+\infty}}^{\phi_{-\infty}}
		(\phi_{-\infty}-z)(z-\phi_{+\infty})
		v_z^2\,dz,
		\]
		and
		\[
		\epsilon\|w_x(t)\|_{L^2(\R)}^2
		\le
		\lambda_4
		\int_{\phi_{+\infty}}^{\phi_{-\infty}}
		(\phi_{-\infty}-z)(z-\phi_{+\infty})
		v_z^2\,dz.
		\]
		Therefore,
		\[
		\frac12\frac{d}{dt}\|w(t)\|_{L^2(\R)}^2
		+
		A_2^*\epsilon\|w_x(t)\|_{L^2(\R)}^2
		\le0,
		\]
		where
		$
		A_2^*
		=
		\frac{2\lambda_3-B_\gamma}{2\lambda_4}.
		$
		Since
		$
		\lambda_3=\frac{\kappa_1m_4}{4},
		$
		the choice
		$
		0<\gamma<
		\sqrt{2(\kappa_1m_4-U_*^2)}-2U_*
		$
		implies
		$
		2\lambda_3-B_\gamma>0,
		$
		and hence
		$
		A_2^*>0.
		$
		Integrating in time yields
		\[
		\|w(t)\|_{L^2(\R)}^2
		+
		2A_2^*\epsilon
		\int_0^t\|w_x(s)\|_{L^2(\R)}^2\,ds
		\le
		\|w(0)\|_{L^2(\R)}^2.
		\]
		
		The above conditional $L^2$-contraction estimate gives
		\[
		\sup_{0\le t\le T}\|w(t)\|_{L^2(\R)}
		\le \|w(0)\|_{L^2(\R)}.
		\]
	Moreover, since $\widetilde{\phi}'\in L^1(\mathbb{R})\cap L^2(\mathbb{R})$ and the bootstrap assumption
	\[
	\sup_{0\le t\le T}\|w(t)\|_{L^\infty(\R)}\le \gamma
	\]
	holds, Lemma~\ref{lem:derivative-estimate} yields
		\[
		\frac{d}{dt}\|w_x(t)\|_{L^2(\R)}^2
		\le
		C_\gamma
		\left(
		\|w_x(t)\|_{L^2(\R)}
		+
		\|w_x(t)\|_{L^2(\R)}^2
		\right),
		\]
		where $C_\gamma>0$ is independent of $T$.
		Lemma~\ref{lem:bootstrap-closure} can therefore
		be applied with $\Gamma=\gamma$ and $c_0=2A_2^*$.
		For sufficiently small
		$
		\|\varphi-\widetilde{\phi}\|_{H^1(\mathbb{R})},
		$
	 the $L^2$-contraction estimate holds for every $t\ge0$, namely,
		\[
		\|w(t)\|_{L^2(\R)}^2
		+
		2A_2^*\epsilon
		\int_0^t\|w_x(s)\|_{L^2(\R)}^2\,ds
		\le
		\|w(0)\|_{L^2(\R)}^2.
		\]
		
\end{proof}
	
	\begin{corollary}[Explicit positive-state version of Theorem~\ref{the:p'=4}]
		\label{cor:explicit-p4}
		Let $\epsilon,\delta>0$ and
		$
		\phi_{-\infty}>\phi_{+\infty}>0.
		$
		Assume that
		\[
		\phi_{+\infty}>y_1,
		\qquad
		0<
		\frac{\delta(\phi_{-\infty}^3-\sigma)}{\epsilon^2}
		\le \frac14,
		\]
		where
		\[
		y_1
		:=
		\frac{
			-8(\sqrt2-1)\phi_{-\infty}
			+
			\sqrt{
				16(\sqrt2-1)[(26-8\sqrt2)\phi_{-\infty}^2
				+12\phi_{-\infty}+3]
			}
		}
		{24(\sqrt2-1)}.
		\]
		Let $\widetilde{\phi}$ be the corresponding strictly decreasing
		shock profile and let $\phi$ be the solution of Eq.~\eqref{eq} with initial
		data $\varphi$. Then there exists $\varepsilon_0>0$ such that, if
		$
		\|\varphi-\widetilde{\phi}\|_{H^1(\mathbb R)}
		<\varepsilon_0,
		$
		the shift $Y(t)$ defined by $Y(0)=0$ and
		\[
		\dot Y(t)
		=
		\frac{c_3}{\phi_{-\infty}-\phi_{+\infty}}
		\int_{\mathbb R}
		\bigl(
		\phi(t,x+Y(t))
		-\widetilde{\phi}(x-\sigma t)
		\bigr)
		\widetilde{\phi}'(x-\sigma t)\,dx,
		\]
		where
		$
		c_3
		=
		-\frac{
			1+4\phi_{-\infty}+6\phi_{-\infty}^2
		}{4}
		$, satisfies, for any $t\in \mathbb{R}_+$,
		\[
		\begin{aligned}
		&	\|\phi(t,\cdot)
			-\widetilde{\phi}(\cdot-\sigma t-Y(t))\|_{L^2(\R)}^2
			+2A_3^*\epsilon
			\int_0^t
			\left\|
			\partial_x
			\bigl[
			\phi(s,\cdot)
			-\widetilde{\phi}(\cdot-\sigma s-Y(s))
			\bigr]
			\right\|_{L^2(\R)}^2\,ds
			\\
			&\le
			\|\varphi-\widetilde{\phi}\|_{L^2(\R)}^2,
		\end{aligned}
		\]
		where
		$
		A_3^*
		=
		\frac{
			8\lambda_3'-
			(1+4\phi_{-\infty}+6\phi_{-\infty}^2)
		}{8\lambda_4'}
		>0.
		$ 
	\end{corollary}
	
	\begin{proof}
		Under the assumption $\phi_{+\infty}>0$, the explicit profile
		estimate \eqref{eq:p=4-inequality} in Corollary~\ref{p=3-4,corollary} applies. The condition
		$\phi_{+\infty}>y_1$ guarantees $A_3^*>0$.
		The conclusion then follows from the proof of Theorem~\ref{the:p'=4} with
		\[
		\gamma=1,
		\qquad
		B_1=
		\frac{1+4\phi_{-\infty}+6\phi_{-\infty}^2}{4},
		\qquad
		\Gamma=1,
		\qquad
		c_0=2A_3^*.
		\]
	\end{proof}

	The preceding contraction estimates provide uniform control of the	perturbation relative to a shifted shock profile. To study the long-time behavior more precisely, we also need quantitative information on the	spatial decay of the profile itself. The preceding two-sided bounds on the profile
	derivative imply that the shock profile
	approaches its asymptotic end states at an
	exponential rate. This result is stated
	precisely in the theorem below.
	
	\begin{theorem}[Exponential decay of the shock profile]
		\label{thm:exponential-decay}
		Let $\widetilde \phi$ be a strictly decreasing shock profile connecting
		$\phi_{-\infty}$ to $\phi_{+\infty}$, with
		$
		\phi_{-\infty}>\phi_{+\infty}$. 
		Assume that there exist constants
		$
		0<\lambda_i\le \lambda_{i+1}
		$
		such that
		\begin{equation}\label{t2}
			\lambda_i
			(\phi_{-\infty}-\widetilde \phi)(\widetilde \phi-\phi_{+\infty})
			\le
			-\epsilon \widetilde \phi'
			\le
			\lambda_{i+1}
			(\phi_{-\infty}-\widetilde \phi)(\widetilde \phi-\phi_{+\infty})
			\qquad \text{on } \mathbb{R}.
		\end{equation}
		Fix any $\xi_0\in\mathbb{R}$ and define
		\[
		s_L:=\widetilde \phi(\xi_0)-\phi_{+\infty}>0,
		\qquad
		s_R:=\phi_{-\infty}-\widetilde \phi(\xi_0)>0,
		\qquad
			D:=s_L+s_R=\phi_{-\infty}-\phi_{+\infty}.
		\]
		Then the following estimates hold:
		for $\xi\le\xi_0$,
		\begin{equation}\label{left}
			s_R
			\exp\left(
			-\frac{\lambda_{i+1}D}{\epsilon}
			|\xi-\xi_0|
			\right)
			\le
			\phi_{-\infty}-\widetilde \phi(\xi)
			\le
			s_R
			\exp\left(
			-\frac{\lambda_i s_L}{\epsilon}
			|\xi-\xi_0|
			\right).
		\end{equation}
		For $\xi\ge\xi_0$,
		\begin{equation}\label{right}
			s_L
			\exp\left(
			-\frac{\lambda_{i+1}D}{\epsilon}
			|\xi-\xi_0|
			\right)
			\le
			\widetilde \phi(\xi)-\phi_{+\infty}
			\le
			s_L
			\exp\left(
			-\frac{\lambda_i s_R}{\epsilon}
			|\xi-\xi_0|
			\right).
		\end{equation}
		In particular, for every $\xi\in\mathbb{R}$,
		\begin{equation}\label{all}
			\lambda_i s_Ls_R
			\exp\left(
			-\frac{\lambda_{i+1}D}{\epsilon}
			|\xi-\xi_0|
			\right)
			\le
			-\epsilon\widetilde \phi'(\xi)
			\le
			\lambda_{i+1}Ds^*
			\exp\left(
			-\frac{\lambda_i s_*}{\epsilon}
			|\xi-\xi_0|
			\right),
		\end{equation}
		where 	\[
		s_*:=\min\{s_L,s_R\},
		\qquad
		s^*:=\max\{s_L,s_R\}.
		\]
	\end{theorem}
	
	\begin{proof}
		Set
		$
		A(\xi):=\phi_{-\infty}-\tilde{\phi}(\xi)$ and 
		$
		B(\xi):=\tilde{\phi}(\xi)-\phi_{+\infty}.
		$
		Then
		$
		A(\xi)+B(\xi)
		=
		\phi_{-\infty}-\phi_{+\infty}
		=
		D.
		$
		We first consider $\xi\le \xi_0$. Since $\tilde{\phi}$ is strictly
		decreasing, $
		\tilde{\phi}(\xi)\ge \tilde{\phi}(\xi_0),
		$
		and hence
		\[
		s_L\le B(\xi)\le D.
		\]
		Moreover,
		$
		A'(\xi)=-\tilde{\phi}'(\xi).
		$
		Therefore, by the two-sided estimate \eqref{t2},
		\[
		\lambda_i A(\xi)B(\xi)
		\le
		\epsilon A'(\xi)
		\le
		\lambda_{i+1}A(\xi)B(\xi),
		\]
		which implies
		\[
		\frac{\lambda_i s_L}{\epsilon}
		\le
		\frac{A'(\xi)}{A(\xi)}
		\le
		\frac{\lambda_{i+1}D}{\epsilon}.
		\]
		Integrating from $\xi$ to $\xi_0$ and using
		$
		A(\xi_0)
		=
		\phi_{-\infty}-\tilde{\phi}(\xi_0)
		=
		s_R,
		$
		one obtains
		\[
		s_R
		\exp\left(
		-\frac{\lambda_{i+1}D}{\epsilon}
		|\xi-\xi_0|
		\right)
		\le
		\phi_{-\infty}-\tilde{\phi}(\xi)
		\le
		s_R
		\exp\left(
		-\frac{\lambda_i s_L}{\epsilon}
		|\xi-\xi_0|
		\right),
		\qquad
		\xi\le\xi_0.
		\]
		This proves \eqref{left}.
		
		For $\xi\ge\xi_0$, monotonicity gives
		$
		s_R\le A(\xi)\le D.
		$
		Since
		$
		B'(\xi)=\tilde{\phi}'(\xi),
		$
		the estimate \eqref{t2} yields
		\[
		\lambda_i A(\xi)B(\xi)
		\le
		-\epsilon B'(\xi)
		\le
		\lambda_{i+1}A(\xi)B(\xi).
		\]
		Consequently,
		\[
		-\frac{\lambda_{i+1}D}{\epsilon}
		\le
		\frac{B'(\xi)}{B(\xi)}
		\le
		-\frac{\lambda_i s_R}{\epsilon}.
		\]
	Integrating from $\xi_0$ to $\xi$ and using
	$B(\xi_0)=\widetilde{\phi}(\xi_0)-\phi_{+\infty}=s_L$ gives
		\[
		s_L
		\exp\left(
		-\frac{\lambda_{i+1}D}{\epsilon}
		|\xi-\xi_0|
		\right)
		\le
		\tilde{\phi}(\xi)-\phi_{+\infty}
		\le
		s_L
		\exp\left(
		-\frac{\lambda_i s_R}{\epsilon}
		|\xi-\xi_0|
		\right),
		\qquad
		\xi\ge\xi_0,
		\]
		which proves \eqref{right}.
		
		It remains to estimate $\tilde{\phi}'$. For $\xi\le\xi_0$,
		using $s_L\le B(\xi)\le D$ together with \eqref{left} and
		\eqref{t2}, it follows that
		\[
		\lambda_i s_Ls_R
		\exp\left(
		-\frac{\lambda_{i+1}D}{\epsilon}
		|\xi-\xi_0|
		\right)
		\le
		-\epsilon\tilde{\phi}'(\xi)
		\le
		\lambda_{i+1}Ds_R
		\exp\left(
		-\frac{\lambda_i s_L}{\epsilon}
		|\xi-\xi_0|
		\right).
		\]
		Similarly, for $\xi\ge\xi_0$, \eqref{right} gives
		\[
		\lambda_i s_Ls_R
		\exp\left(
		-\frac{\lambda_{i+1}D}{\epsilon}
		|\xi-\xi_0|
		\right)
		\le
		-\epsilon\tilde{\phi}'(\xi)
		\le
		\lambda_{i+1}Ds_L
		\exp\left(
		-\frac{\lambda_i s_R}{\epsilon}
		|\xi-\xi_0|
		\right).
		\]
		Combining the preceding estimates, we obtain for every
		$\xi\in\mathbb R$,
		\[
		\lambda_i s_Ls_R
		\exp\left(
		-\frac{\lambda_{i+1}D}{\epsilon}
		|\xi-\xi_0|
		\right)
		\le
		-\epsilon\tilde{\phi}'(\xi)
		\le
		\lambda_{i+1}Ds^*
		\exp\left(
		-\frac{\lambda_i s_*}{\epsilon}
		|\xi-\xi_0|
		\right),
		\]
		which is \eqref{all}.
	\end{proof}

	The exponential decay of the shock profile established in Theorem~\ref{thm:exponential-decay}
	implies, in particular, that
	\[
	\widetilde{\phi}'\in L^1(\mathbb{R})\cap L^2(\mathbb{R}).
	\]
	Combined with the global $L^2$-contraction estimates obtained above,
	this allows us to control the spatial derivative of the perturbation and
	derive its decay as $t\to\infty$. We first establish this fact in the
	following lemma, and then deduce the time-asymptotic stability and the
	decay of the shift.
	
	\begin{lemma}[Decay of the spatial derivative]\label{lem:spatial-derivative-decay}
		For either $p=3$ or $p=4$, consider the
		global solution $\phi$ of \eqref{eq}
		and its associated strictly decreasing
		shock profile $\widetilde{\phi}$.
		Suppose that the assumptions specified in
		Theorem~\ref{thm:main} hold in the cubic case,
		while those in Theorem~\ref{the:p'=4}
		are fulfilled in the quartic case.
		Denote by $Y(t)$ the shift function
		introduced in the corresponding theorem. Then
		\[
		\left\|
		\partial_x\left[
		\phi(t,\cdot)
		-\widetilde{\phi}(\cdot-\sigma t-Y(t))
		\right]
		\right\|_{L^2(\mathbb R)}
		\to0,
		\qquad t\to\infty.
		\]
	\end{lemma}
	
	\begin{proof}
		Set
		$
		w(t,x)
		:=
		\phi(t,x)-\widetilde{\phi}(x-\sigma t-Y(t))$ and 
		$
		D(t):=\|w_x(t)\|_{L^2(\R)}^2.
		$
		The corresponding $L^2$-contraction estimate gives
		\[
		\sup_{t\ge0}\|w(t)\|_{L^2(\R)}<\infty,
		\qquad
		\int_0^\infty D(t)\,dt<\infty.
		\]
		Moreover, the assumptions of the corresponding contraction theorem
		imply
		\[
		\sup_{t\ge0}\|w(t)\|_{L^\infty(\R)}<\infty.
		\]
		Since $\widetilde{\phi'}\in L^1(\mathbb{R})\cap L^2(\mathbb{R}) $, the definition of the shift also yields
		\[
		|\dot Y(t)|
		\le C\|w(t)\|_{L^\infty(\R)}.
		\]
		Hence Lemma~\ref{lem:derivative-estimate} implies
		\[
		D'(t)
		\le
		C\left(D(t)^{1/2}+D(t)\right),
		\qquad t\ge0.
		\]
		
		It remains to prove that
		\[
		D(t)\to0,
		\qquad t\to\infty.
		\]
		Suppose otherwise. Since $D$ is continuous and
		$D\in L^1(0,\infty)$, there exist $\eta>0$ and infinitely many
		pairwise disjoint intervals $[s_n,t_n]$ such that
		$
		D(s_n)=\frac{\eta}{2}$,
		$
		D(t_n)=\eta,
		$
		and
		\[
		\frac{\eta}{2}\le D(t)\le\eta,
		\qquad t\in[s_n,t_n].
		\]
		Thus,
		\[
		D'(t)
		\le C\left(\eta^{1/2}+\eta\right)
		\qquad\text{on }[s_n,t_n],
		\]
		and after integrating both sides, it follows that
		\[
		t_n-s_n
		\ge
		\frac{\eta}
		{2C(\eta^{1/2}+\eta)}
		=:c_\eta>0.
		\]
		Therefore,
		\[
		\int_{s_n}^{t_n}D(t)\,dt
		\ge
		\frac{\eta}{2}c_\eta>0.
		\]
		Since the intervals $[s_n,t_n]$ are pairwise disjoint, this
		contradicts
		\[
		\int_0^\infty D(t)\,dt<\infty.
		\]
		Consequently,
		$
		D(t)\to0$ ( 
		$ t\to\infty
		$), 
		and therefore
		\[
		\|w_x(t)\|_{L^2(\mathbb{R})}
		\to0,
		\qquad t\to\infty.
		\]
	\end{proof}

\begin{proof}	
({\bf Proof of Theorem \ref{thm:asymptotic-stability}})
		Set
		$
		w(t,x):=
		\phi(t,x)
		-\widetilde{\phi}(x-\sigma t-Y(t)).
		$
		By Lemma~\ref{lem:spatial-derivative-decay},
		\[
		\|w_x(t)\|_{L^2(\R)}\to0,
		\qquad t\to\infty.
		\]
		The $L^2$-contraction estimate in Theorem \ref{thm:main} or Theorem \ref{the:p'=4} gives
		\[
		\sup_{t\geq 0}\|w(t)\|_{L^2(\R)}<\infty.
		\]
		Therefore, by the one-dimensional Gagliardo-Nirenberg inequality,
		\[
		\|w(t)\|_{L^\infty(\R)}^2
		\leq
		2\|w(t)\|_{L^2(\R)}\|w_x(t)\|_{L^2(\R)}
		\to0.
		\]
		For $2<q<\infty$, interpolation yields
		\[
		\|w(t)\|_{L^q}
		\leq
		\|w(t)\|_{L^\infty(\R)}^{1-\frac{2}{q}}
		\|w(t)\|_{L^2(\R)}^{\frac{2}{q}},
		\]
		and we infer that
		\[
		\|w(t)\|_{L^q}\to0,
		\qquad 2<q\leq\infty.
		\]
		Next, by the definition of the corresponding shift,
		\[
		\begin{aligned}
		|\dot{Y}(t)|
		&=
		\left| \frac{c_i}{\phi_{-\infty}-\phi_{+\infty}}
		\int_{\mathbb{R}} w(t,x)\,\phi'\bigl(x-\sigma t-Y(t)\bigr)\,dx\right| \\
		&\le
		 \frac{|c_i|}{\phi_{-\infty}-\phi_{+\infty}}
		\|w(t)\|_{L^\infty(\mathbb{R})}
		\|\phi'\|_{L^1(\mathbb{R})}, 
	\end{aligned}
		\]
		where \(c_i=c_1\) when \(p=3\), and \(c_i=c_2\) when \(p=4\). 
		Since $\widetilde{\phi}'\in L^1(\mathbb{R})$ and
		$\|w(t)\|_{L^\infty(\R)}\to0$, we conclude that
		\[
		|\dot Y(t)|\to0,
		\qquad t\to\infty.
		\]
	Finally, since \(Y(0)=0\), it follows that
	\[
	\frac{Y(t)}{t}
	=
	\frac{1}{t}\int_0^t \dot{Y}(s)\,ds.
	\]
	For any \(\varepsilon>0\), the convergence \(\dot{Y}(t)\to 0\) implies that there exists \(T_\varepsilon>0\) such that
	\[
	|\dot{Y}(s)|<\varepsilon,\qquad s\ge T_\varepsilon.
	\]
	Hence, for \(t>T_\varepsilon\),
	\[
	\left|\frac{Y(t)}{t}\right|
	\le \frac{1}{t}\int_0^{T_\varepsilon}|\dot{Y}(s)|\,ds
	+\frac{1}{t}\int_{T_\varepsilon}^t|\dot{Y}(s)|\,ds
	\le \frac{1}{t}\int_0^{T_\varepsilon}|\dot{Y}(s)|\,ds+\varepsilon.
	\]
	Letting \(t\to\infty\) and then \(\varepsilon\to0\), we obtain
	\[
	\frac{Y(t)}{t}\to 0,\qquad t\to\infty.
	\]

\end{proof}






\bibliographystyle{plain}

\begin{thebibliography}{10}

\bibitem{Boussinesq} Boussinesq, J.: Th\'{e}orie des ondes et des remous qui se propagent le long d'un canal rectangulaire horizontal, en communiquant au liquide contenu dans ce canal des vitesses sensiblement pareilles de la surface au fond. \emph{J. Math. Pures Appl.}, 17, 55-108 (1872). 

\bibitem{DJkorteweg1895} Korteweg, D. J. and de Vries, G.: XLI. On the change of form of long waves advancing in a rectangular canal, and on a new type of long stationary waves. \emph{Philos. Mag.}, 39(240), 422-443 (1895).

\bibitem{SuGardner1969} Su, C. H. and Gardner, C. S.: Korteweg–de Vries equation and generalizations. III. Derivation of the Korteweg–de Vries equation and Burgers equation. \emph{J. Math. Phys.}, 10, 536-539 (1969).

\bibitem{p3Demiray1999} Demiray, H.: Propagation of weakly nonlinear waves in fluid-filled thick viscoelastic tubes. \emph{Appl. Math. Model.} 23(10), 779-798 (1999). \DOI{10.1016/S0307-904X(99)00012-8.}

\bibitem{p3Verheest2019} Verheest, F. and Hereman, W. A.: Collisions of acoustic solitons and their electric fields in plasmas at critical compositions, \emph{J. Plasma Phys.} 85(1), 905850106 (2019). \DOI{10.1017/S0022377818001368}

\bibitem{p3PUXI2021} Pu, X. and Xi, X.: Derivation of the mKdV equation from the Euler–Poisson system at critical densities. \emph{J. Differential Equations}, 282, 446-480 (2021). \DOI{10.1016/j.jde.2021.02.026.}



\bibitem{p3Sarma2009} Sarma, A.: Ion-acoustic shock waves with negative ions in presence of dust particulates. \emph{Phys. Lett. A}, 373(45), 4174-4177 (2009). \DOI{10.1016/j.physleta.2009.09.034.}

\bibitem{p4Verheest2016} Verheest, F., Olivier, C. P. and Hereman, W. A.: Modified Korteweg–de Vries solitons at supercritical densities in two-electron temperature plasmas. \emph{J. Plasma Phys}, 82(2) 905820208 (2016). \DOI{10.1017/S0022377816000349.}


\bibitem{Molinetlowregularity2002} Molinet, F. and Ribaud, F.: On the low regularity of the Korteweg–de Vries-Burgers equation. \emph{Int. Math. Res. Not.}, 2002(37), 1979-2005 (2002). \DOI{10.1155/S1073792802112104.}


\bibitem{Zhang2009} Zhang, H. and Han, L.: Global well-posedness and inviscid limit for the modified Korteweg–de Vries-Burgers equation. \emph{Nonlinear Anal.}, 71(12), e1708-e1715 (2009). \DOI{10.1016/j.na.2009.02.042.}

\bibitem{Gallego2019} Gallego, F. A. and Pazoto, A. F.: On the well-posedness and asymptotic behaviour of the generalized Korteweg–de Vries-Burgers equation. \emph{Proc. Roy. Soc. Edinburgh Sect. A}, 149(1), 219-260 (2019). \DOI{10.1017/S0308210518000240.}


\bibitem{chen2026}Chen, G., Eun, N., Kang, M.-J. and Shen, Y.: Uniform stability of oscillatory shocks for KdV-Burgers equation. \emph{arXiv preprint}, arXiv:2602.22652 (2026).

\bibitem{Chen2000}
Chen, H.: 
Travelling-waves for generalized KdV-Burgers equations. 
In: Bona, J. L., Saxton, K., Saxton, R. (eds.) 
\emph{Nonlinear PDE's, Dynamics and Continuum Physics}, 
Contemporary Mathematics, vol. 255, pp. 1-24. 
American Mathematical Society, Providence, RI (2000).

\bibitem{kang2017}Kang, M.-J. and Vasseur, A. F.: $L^2$-contraction for shock waves of scalar viscous conservation laws. \emph{Ann. Inst. H. Poincaré C Anal. Non Linéaire}, 34(1), 139-156 (2017). \DOI{10.1016/j.anihpc.2015.10.004.}

\bibitem{chen20262}Chen, G., Eun, N., Kang, M.-J. and Shen, Y.: $L^2$-contraction of shock waves for KdV-Burgers equation. \emph{arXiv preprint}, arXiv:2603.09263 (2026). 


\bibitem{Bona1985}Bona, J. L. and Schonbek, M. E.: Travelling-wave solutions to the Korteweg-de Vries-Burgers equation. \emph{Proc. Roy. Soc. Edinburgh Sect. A}, 101(3-4), 207-226 (1985). \DOI{10.1017/S0308210500020783.} 

\bibitem{Pego1985}Pego, R. L.: Remarks on the stability of shock profiles for conservation laws with dissipation. \emph{Trans. Amer. Math. Soc.}, 291(1), 353-361 (1985). \DOI{10.1090/S0002-9947-1985-0797065-0.}

\bibitem{Pego1993}Pego, R. L., Smereka, P., and Weinstein, M. I.: Oscillatory instability of traveling waves for a KdV-Burgers equation. \emph{Phys. D}, 67(1-3), 45-65 (1993). \DOI{10.1016/0167-2789(93)90197-9.}

\bibitem{zhang2019}Zhang, W.-G., Li, W.-X., Deng, S.-E. and Li, X.: Asymptotic stability of monotone decreasing kink profile solitary wave solutions for generalized KdV-Burgers equation. \emph{Acta Math. Appl. Sin. Engl. Ser.}, 35, 475-490 (2019). \DOI{10.1007/s10255-019-0825-7.}

\bibitem{Eun2026} Eun, N., Han, S. and Kim, J.: Large-time behavior towards composite waves of degenerate shock and rarefaction wave for modified KdV-Burgers equation. \emph{arXiv preprint}, arXiv:2609.20591 (2026).





\bibitem{MR4195742}Kang, M.-J. and Vasseur, A. F.: Contraction property for large perturbations of shocks of the barotropic Navier-Stokes system. \emph{J. Eur. Math. Soc. (JEMS)}, 23(2), 585-638 (2021). \DOI{10.4171/jems/1018.}



\end{thebibliography}

 \end{document}